%% file: paperForum.tex
\input artmacros

\advance\vsize by 5 mm
\advance\voffset by -5 mm
%
%% 
%%%%%%%%%%%%%% SIMBOLI
%
\mathdef\G{\mathchar"100}%
\mathdef\a{\alpha}%
\def\aut#1{\mathord{\mathop {\rm Aut}#1}}%
\mathdef\autg{\aut G}%
\def\caut#1{\mathord{\mathop {\rm CAut}#1}}%
\mathdef\cautg{\caut G}%
\def\End#1{\mathord{\mathop {\rm End}#1}}%
\def\CEnd{\mathop {\rm CEnd}}
\mathdef\CEndg{\CEnd G}% 
\def\inn#1{\mathord{\mathop {\rm Inn}#1}}%
\def\tor#1{\mathord{\mathop {\rm tor}#1}}%
\def\Hom{\mathop{\rm Hom}\nolimits}%
\def\Der{\mathop{\rm Der}\nolimits}%
\def\C{{\tondo C}}%
\def\pru{{\C_{p^\infty}}}%
\def\Dr{\mathop{\rm Dr}}%
\def\Cr{\mathop{\rm Cr}}%
\mathdef\Sy{{\Tondo S}}%
\def\smatrix(#1,#2,#3,#4){\bigl({#1\atop #3}{#2\atop #4}\)}% % %
\relchardef\shepi=\amsafam16%
%%
%\font\titfont=pplr8r at 16 pt
\font\titfont=cmr12 at 18 pt
\newfam\frakfam
\font\tenfrak=eufm10
\textfont\frakfam\tenfrak
\def\frak{\tenfrak\fam\frakfam\relax}
\def\restr#1{_{ |#1}}%
% 
%	BIBLIOGRAFIA
%
\ListOfReferences 
[B]	% Baer
%[BW]	% Buckley-Wiegold
[D]	% Dixon
[DG]	% Dugas-G"obel
[Ku] 	% Kueker
[EN]	% Hall, Edmonton Notes
[H]	% Hulse
[OP]	% Puglisi
[RobinH]	% Robinson
[RobinStA]	% Robinson, St. Andrews
[Rot]	% Rotman, libro
[ST]	% Thomas
%[VW]	% Vaughan-Lee-Wiegold
[Z]	% Zaleskii
\EndOfReferences
\def\[#1]{\bibref{#1}}%
\tracingpages=2 %
\centerline{\titfont Subgroups defining automorphisms in} \medskip
\centerline{\titfont locally nilpotent groups} \vskip 1 cm
\centerline{\caps Giovanni Cutolo and Chiara Nicotera} \vskip 1 cm
{\narrower\parindent=0pt\caps
Abstract:\enskip\rm We investigate some situation in which 
automorphisms of a group~$G$ are unique\-ly determined by their 
restrictions to a proper subgroup~$H$. Much of the paper is devoted  
to studying under which additional hypotheses this property forces $G$ to
be nilpotent if $H$ is.  As an application we prove that certain
countably infinite locally nilpotent groups have uncountably many
(outer) automorphisms.

\bigskip \caps
\smallskip
\caps Math.~Subj.~Classification (2000):\enskip\rm 20F28, 20F19.\par}
\vskip 1 cm
\noindent
Homomorphisms of
groups are defined by their restriction to any generating set of 
their domain.  This property actually
characterizes generating subsets of groups, for if~$H$ is a proper 
subgroup of the group~$G$ then there are two different homomorphisms 
from~$G$ to the same group~$K$ whose restrictions to~$H$ are the 
same: a simple construction due to Eilenberg and Moore is suggested 
as Exercise~3.35  in~\[Rot], p.\thinspace54.

The situation can be quite different if---rather than referring 
to all homomorphisms of domain a group~$G$---we restrict attention 
to, say, endomorphisms of~$G$ only. For instance, if $G$ is 
isomorphic to  the rational group~$\Q$ and $g$ is any nontrivial 
element of~$G$, then two endomorphisms of~$G$ coincide if they agree 
on~$g$, in other words  endomorphisms of~$G$ are uniquely determined by
their 
restrictions to~$\{g\}$. This suggests the following definition. Let 
$G$ be a group and let~\G be  a set of endomorphisms of~$G$.
We say that a subset~$X$ of~$G$ is a {\it $\G$\@basis} if and only 
if,  for every $\a,\beta\in\G$, it holds $\a=\beta$ if $\a\restr 
X=\beta\restr  X$, where $\a\restr X$ and $\beta\restr X$ denote 
restrictions to~$X$. We shall also use expressions like `End-basis', 
`Aut-basis' or `Inn-basis of~$G$' as synonym with $\End G$\@, 
$\autg$\@, or $\inn G$\@basis  respectively.
For instance, the above example can be rephrased by saying that 
in the rational group every one-element subset 
different from the identity subgroup is an End-basis. More  generally,
every maximal independent subset of a torsion-free abelian  group~$A$ 
is an End-basis of~$A$.

The Aut-bases of a group~$G$ are just the bases of~$\autg$ viewed as a
permutation group on~$G$, whence the name.  Indeed, it is clear that
for any $\G\le \autg$ a subset $X$ of~$G$ is a $\G$\@basis if and only
if $C_\G(X)=1$.  In particular, the Inn-bases of a group~$G$ are the
subsets~$X$ of~$G$ such that $C_G(X)=Z(G)$.

Other self-evident facts about $\G$\@bases (for a set~\G of 
endomorphisms of a group~$G$) are that if $X$ is a $\G$\@basis then $X_1$ 
is a $\G_1$\@basis for any subset $\G_1$ of~\G and any superset $X_1$ 
of~$X$ contained in~$G$.  Also, $X$ is a $\G$\@basis if and only if 
$\gen X$ is a $\G$\@basis, so the property of being a $\G$\@basis 
could be equivalently defined as an embedding property for subgroups. 
In this paper we shall assume this point of view and discuss some  
instances of the general problem 
of when a group~$G$ inherits group theoretical properties from  a 
subgroup of~$G$ which is a $\G$\@basis, for some specific $\G\subseteq 
\End G$. For example, it is immediate to observe that a group is 
abelian if it has an abelian subgroup which is an Inn-basis (Lemma 
\rif{abel}).
We will be mainly concerned with the case $\G=\autg$.

A drastically restrictive result is that a direct power of
every centreless group can be 
embedded as a normal subgroup which is an Aut-basis in a group with 
rather arbitrary structure (see Corollary\rif{costruzione}).  This is 
the reason why we turn our attention to group classes without 
nontrivial centreless groups, and mainly study nilpotent (sub)groups. 

Even the property of having a nilpotent subgroup as an  Aut-basis 
seems not to be very strong by itself. We give several examples of 
non-nilpotent groups with a nilpotent subgroup as an Aut-basis; such 
groups may even be locally nilpotent and the  subgroup may satisfy 
various embedding conditions.
However, it is possible to obtain some information in the positive.
The main results in the first two sections of this paper are the following.
Assume that $H$ is a subgroup and an Aut-basis of the group~$G$.
If $G$ is nilpotent then it has the same nilpotency class as~$H$
(see Theorem\rif{Loryno}).  A sufficient condition for $G$ to be nilpotent is that
$H$ is a subnormal nilpotent subgroup of the hypercentre of~$G$
(Theorem\rif{Loryno+}).  Another sufficient condition is given in
Theorem\rif{Teorema Pasqualino}: if $G$ is locally nilpotent and
$H$ is not only an Aut-basis of~$G$ but also an Aut-basis of every
subgroup of~$G$ containing it, then $G$ is hypercentral, or nilpotent,
provided $H$ has the same property. In this case, also, the 
hypercentral lengths of $H$ and~$G$ coincide.

A condition of special interest to us is when the basis considered is
finite, or, equivalently, the subgroup considered is finitely
generated.  Results about the existence of finite Aut-bases in groups
(or other structures) have appeared in the literature.  By a theorem in
model theory due to D.W.~Kueker~\[Ku] a countable group~$G$ has a
finite Aut-basis if and only if \autg is countable, and if this fails
to happen then $|\autg|=2^{\aleph_0}$ (a proof only using the language
of group theory is in~\[B], Lemma~III.1).  
% We mention in passing that
% obvious adaptations of the proof yield the same result---a countable
% group~$G$ has a finite $\G$\@basis if and only if \G is countable,
% otherwise $|\G|=2^{\aleph_0}$---for each of the following choices
% for~$\G$ too: $\End G$, the set of monic endomorphisms of~$G$, the set
% of epic endomorphisms of~$G$. 
A related result by R.~Baer is also in~\[B], Satz~III.4: every quotient of a
group~$G$ has countable automorphism group if and only if
$G$ is countable and every quotient of~$G$ has a finite Aut-basis. 
Relevant examples of groups with a finite Aut-basis are the groups
with only finitely many automorphisms. Such groups have been the subject
of several investigations (see for instance~\[RobinH]).  Among them
there exist torsion-free abelian groups of very high cardinality (more
precisely, of any arbitrary infinite cardinality less than the first
strongly inaccessible cardinal) whose only endomorphisms are the
mappings $x\mapsto x^n$ for integers~$n$, so that the only non-trivial
automorphism is the inverting automorphism.  Of course the singleton
of any non-trivial element of such a group is an End-basis of it.

Locally nilpotent groups with a finite Aut-basis are the subject of
the third and last section of this paper.  A special case of
Theorem\rif{finito} is that a locally nilpotent group has a finite
subgroup as an Aut-basis if and only if it is finite.  From the same
Theorem\rif{finito} we also draw some consequences on the existence of
outer automorphisms.  It is a theorem by O.~Puglisi~\[OP] that every
periodic countably infinite locally nilpotent group has~$2^{\aleph_0}$
(outer) automorphisms.  We extend this result to countable
periodic-by-(finitely generated) locally nilpotent groups which are
not finitely generated (Corollary\rif{outer}).

\Sezione First results and examples

For convenience of further reference we start by stating again an  
observation already made in the introduction.

\Lemma\detto{centro} A subset~$X$ of a group~$G$ is an Inn-basis if 
and only if $C_G(X)=Z(G)$. In particular, the latter equality  holds 
if $X$ is an Aut-basis.

\Lemma\detto{condHom} Let the subgroup $H$ be an Aut-basis of the
group~$G$ and let
$N=H^G$. Then $\Hom\(G/N, Z(H)\)=0$.

\Dim 
The group of all automorphisms of~$G$ acting  trivially on both~$N$ 
and~$G/N$ is trivial since $N$ is an Aut-basis of~$G$. As is 
well-known (see \[RobinStA], p.\thinspace 66, for instance),
this group is isomorphic to the group  $\Der\(G/N, Z(N)\)$ 
of all derivations from $G/N$ to~$Z(N)$. Now, Lemma\rif{centro} 
yields $\Der\(G/N, Z(N)\)=\Hom\(G/N, Z(N)\)$, so the  latter group is 
trivial as well. Finally $Z(H)=Z(G)\cap H\le Z(N)$ by Lemma\rif{centro},
hence $\Hom\(G/N, Z(H)\)=0$. \qed

A special case of the next lemma is the well-known fact
that the property of being an 
Aut-basis and the property of being an 
Inn-basis are equivalent for a normal subset of a centreless group.

\Lemma\detto{AutBasiNormali}
Let $N$ be a normal subgroup of a group~$G$ containing $Z(G)$.  Then 
$N$ is an Aut-basis of~$G$ if and only if  $C_G(N)=Z(G)$ and  
$\Hom\(G/N,Z(G)\)=0$.

\Dim 
The two previous lemmas ensure that the condition is necessary. 
Conversely, assuming that it holds, we have to check that  
$\G:=C_{\autg}(N)$ is trivial. By the Three Subgroup Lemma (see for 
instance \[EN], Lemma 3.2) $[G,\G,N]=1$, so, by hypothesis,
$[G,\G]\le C_G(N)=Z(G)\le N$. Thus
\G acts trivially on both $N$ and $G/N$. Hence \G embeds in
$\Der\(G/N, Z(N)\)=\Hom\(G/N,Z(G)\)=0$.\qed

Now, if $N\n G$ and $C_G(N)=1$ then $C_H(N)=1=Z(H)$ for all  
subgroups $H$ of~$G$ containing~$N$. Hence $N$ is an Aut-basis of  
every such~$H$. This situation suggests the following definition.

Let us say that a subset $X$ of a group~$G$ is a {\it hereditary  
Aut-basis} of~$G$ if it is an Aut-basis of every subgroup of~$G$  
containing it. (This is still a property that generalizes the  
property of being a generating set.)

Thus Lemma\rif{AutBasiNormali} shows that a normal subgroup of a group~$G$
whose centralizer in~$G$ is trivial is a hereditary Aut-basis. 
This is the statement we will make use of, however we note that from
Lemma~8.1.1 of~\[H] the following stronger result follows:

\Lemma\detto{centroidentico}
Let $H$ be an ascendant subgroup of a group~$G$ such that $C_G(H)=1$.  
Then $H$ is a hereditary Aut-basis of~$G$.

Therefore every centreless group embeds as an ascendant subgroup and
hereditary Aut-basis in a complete group~$\hat H$, by means of the
automorphism tower construction.  If $H$ is finite then~$\hat H$ is also
finite and $H$ is subnormal in~$\hat H$.  Another embedding, of a
group closely related to a given centreless group, as a normal
hereditary Aut-basis in a much more arbitrary group is given by the
following corollary.

\Cor\detto{costruzione}
Let $H$ and $K$ be groups and  assume  $Z(H)=1$. Then there exists
a group~$G$  containing a normal subgroup~$N$ such that:
\\i) $N$ is isomorphic to the direct product of $\kappa$ copies
of~$H$, where $\kappa=\max\{\aleph_0,|K|\}$;
\\ii) $G/N\iso K$;
\\iii) $N$ is a hereditary Aut-basis of~$G$.

\Pf
If $K$ is finite let $K^*$ be a countably infinite group in which
$K$ is embedded, otherwise let $K^*=K$.  Hence $\kappa=|K^*|$.  Now
let $G^*$ be the standard wreath product of~$H$ by~$K^*$ and let $N$ be
its base subgroup.  In view of Lemma\rif{centroidentico} the
required conditions are satisfied if we set $G=NK$.\qed

Now, we are looking for statements of the following kind: ``if the 
group $G$  has a subgroup $H$ which is an Aut-basis and satisfies 
some group theoretical condition~$\frak X$ then $G$ itself  
satisfies~$\frak X$". Possibly we could have to strengthen the  
hypothesis by also requiring some embedding property for~$H$, or~that 
$G$ satisfies some condition weaker than~$\frak X$. 
Corollary\rif{costruzione} suggests that reasonable  
candidates for~$\frak X$ in order that such a statement be true---at  
least under the hypothesis that $H$ is normal in~$G$---are finiteness 
conditions and conditions, like nilpotency or  hypercentrality, which 
exclude nontrivial centreless groups.

Rather than on finiteness conditions our interest will be focused on 
nilpotent subgroups as Aut-bases.
With regard to finiteness conditions we only recall
how well-known results on groups of automorphism of groups 
satisfying finiteness conditions translate into analogous results 
about normal Inn-bases, and observe that normality is necessary here.
Indeed, if the normal subgroup~$N$ of a group~$G$ 
is an Inn-basis then $G/Z(G)=G/C_G(N)$ embeds in $\aut N$.  For 
instance, if $N$ is polycyclic and $G$ is a radical group, then $G$ 
is centre-by-polycyclic. In particular, {\sl every locally  
nilpotent group with a finitely generated normal subgroup as an
Inn-basis is nilpotent}. 
As said, normality of~$N$ plays a relevant role in this context, 
even if the condition of being an Inn-basis is replaced by the
stronger requirement that $N$ is an Aut-basis. Indeed, 
let $A$ be one of the torsion-free abelian 
groups of high cardinality referred to in the introduction,
whose only non-trivial automorphism~\a is 
the inversion, and let $G=A\semid\gen\a$. For any $a\in A\setminus 1$ 
then  $H=\gen{a,\a}$ is an Aut-basis 
of $G$. Clearly $H$ is isomorphic to the infinite dihedral group, 
while $G$ is a metabelian centreless group of the same cardinality 
as~$A$.
\smallskip
The next example shows that a group does not need
to be locally nilpotent, or soluble, even if it has a normal 
nilpotent subgroup which is a finitely generated  Aut-basis.

\Example\detto{ex1}
Let $N$ be the free group on two elements $x$~and $y$ in the 
variety of nilpotent groups of class at most~2. Then $Z(N)=N'$ is 
generated  by~$c=[x,y]$. Let $\gamma$ be
the  involutory automorphism of~$N$ defined by $x\mapsto x^{-1}$ and  
$y\mapsto y^{-1}$, and let $G=N\semid\gen\gamma$. 
Then  $C_G(N)=\gen c=Z(G)$ and $\Hom\(G/N,Z(G)\)=0$. Thus we may apply 
Lemma\rif{AutBasiNormali} to obtain that $N$ is an  
Aut-basis of the non-nilpotent  polycyclic group~$G$. As $N$ is maximal
in~$G$, it is even a hereditary Aut-basis.

An insoluble example can be constructed as follows.
The action of $\aut N$ on $N_{\rm ab}=N/N'$ gives rise to an 
epimorphism $\phhi:\aut N\epi {\rm GL}_2(\Z)$, since every 
automorphism of~$N_{\rm ab}$ is induced by an automorphism of~$N$.
The kernel of~$\phhi$ is $I:=\inn N$. Consider the group~$H=N\semid 
\G$, where $\G$ is the centralizer of~$c$ 
in~$\aut N$. As immediately checked, $\G$
is the preimage of~${\rm SL}_2(\Z)$ under~$\phhi$, thus
$\G/I\iso {\rm SL}_2(\Z)$.
Let $C=\{a^{-1}\tilde a\st a\in N\}$,
where $\tilde a$ denotes the inner 
automorphism of~$N$ determined by~$a$.  Then $C=C_H(N)\iso N$ and 
$C_H(C)=N$. Clearly $NC=NI$ is nilpotent of class~$2$.
Now, $NC$ is an Aut-basis of~$H$: this follows again from 
Lemma\rif{AutBasiNormali}, since $C_H(NC)=C\cap N=\gen c=Z(H)$ and
$\Hom\(H/NC,Z(H)\)\iso\Hom \(\G/I,\gen c\)=0$.
\qed

Despite the above example there are some cases in which the presence of a 
nilpotent  subgroup that is an Aut-basis in a group forces the 
latter to be nilpotent. We first record a simple
remark.

\Lemma\detto{abel} Let the subgroup~$H$ be an  Inn-basis of the 
group~$G$. Then, for any $n\in\N$:
\\i) $G$ is nilpotent of class at 
most~$n$ if and only if $H\le  Z_n(G)$;
\\ii) if $H$ is abelian then $G$ is abelian;
\\iii) if $H\le G'$ and $G$ is hypercentral then $G$ is abelian.

\Pf
Assume $H\le Z_n(G)$. Then $Z(G)=C_G(H)\ge 
C_G\(Z_n(G)\)\ge\gamma_n(G)$, so $\gamma_{n+1}(G)=1$. This 
proves~$(i)$. To prove~$(ii)$ suppose  that $H$ is abelian. Then 
$H\le C_G(H)=Z(G)$ and $(i)$ shows that  $G$ is abelian.
Finally, assume $H\le G'$. Then $Z_2(G)\le C_G(G')=Z(G)$. If $G$ is 
hypercentral this shows that $G$ is abelian.\qed

\Prop\detto{metab}
Let $G$ be a metabelian group and let 
$H$ be a subnormal subgroup of defect~$d$ in~$G$.
If $H$ is an Inn-basis 
of~$G$ then $Z_i(H)\le Z_{i+d}(G)$ for every positive integer~$i$.
In particular,
% $Z_\omega(H)= H\cap Z_\omega(G)$, and 
if $H$ is nilpotent 
of class~$c$ then $G$ is nilpotent of class at most~$c+d$.

\Pf
If $X$ is any metabelian group, then $[Y,a,b]=[Y,b,a]$ for all $a,b\in 
X$ and $Y\le X'$, as the automorphisms induced by conjugation on~$X'$ 
by~$a$~and~$b$ generate a commutative subring of~$\End X'$.
Therefore $[U,V,L_1,L_2,\ldots,L_t]=[U,V,L_{1\sigma},L_{2\sigma},
\ldots,L_{t\sigma}]$ for any $t\in\N$, $\;U,V,L_1,\ldots,L_t\le X$
and $\sigma\in\S_t$.

To prove the proposition it will be clearly enough to assume $d=1$, 
that is:  $H\n G$.
Let $K=Z_i(H)$. Then $[K,G]\le K$ and $[K,G,{}_iH]=1$.
Hence $[K,G,{}_{i-1}H]\le C_G(H)=Z(G)$, so $[K,G,{}_{i-1}H,G]=1$.
If $i>1$ we can apply the property of metabelian groups just recalled
to get $[K,G,{}_{i-2}H,G,H]=1$.
This again gives $[K,G,{}_{i-2}H,G,G]=1$, because $H$ is an Inn-basis.
By carrying on this procedure we eventually obtain $[K,H,{}_iH]=1$.
Therefore $K\le Z_{i+1}(G)$, as required. If $H$ is nilpotent
the last clause in the statement follows now from 
Lemma\rif{abel}$\;(i)$.
\qed

Some remarks about Proposition\rif{metab} are in order. Firstly, every 
nonabelian subgroup of a dihedral group is an Inn-basis. This makes 
easy to find examples showing that nonnilpotent metabelian groups may 
have (non-subnormal) nilpotent subgroups as Inn-bases and that the 
bound $c+d$ for the nilpotency class of~$G$ found in the proposition 
is the best possible (unless $c=1$, of course, in which case $G$ is abelian
by Lemma\rif{abel}). 
Secondly, 
% the statement does not hold for ordinals grater 
% than~$\omega$. Indeed, let $p$ be a prime, let $A$ be an abelian 
% $p$\@group of infinite exponent and let $B$ be a nontrivial abelian group.
% Let $G= (A\times B)\semid \gen x$, where $\gen x$ acts 
% fixed-point-freely on~$B$ and faithfully on~$A$, and $a^x=a^{p+1}$
% for all $a\in A$. Then $N:=A\gen x$ is a normal 
% Inn-basis of~$G$, but $N=Z_{\omega+1}(N)$ is not contained in the 
% hypercentre of~$G$, which is~$A$. Finally,
% we note that 
the hypothesis that the group be metabelian is essential.
% in Proposition\rif{metab}.
It is in fact possible to construct centre-by-metabelian groups
with a nilpotent normal subgroup of class~$2$ as an Inn-basis and which are either
non-nilpotent (see Example\rif{ex1}, or also Examples\rif{27} and\rif{22})
or nilpotent of arbitrarily high class.
Indeed, let $n\in\N$ and let $W$ be the standard wreath product of 
the quaternion group of order~$8$ and a cyclic group of order~$2^n$.  
Let $B$ be the base group of~$W$ and $C=[W,Z(B)]$.  Then $B/C$ is 
an Inn-basis of~$G=W/C$ and has class~$2$, while the class of~$G$ is 
at least~$2^n+1$. Further examples with the same properties can be 
read off from Example\rif{22}, in particular 
they can be chosen to be $p$\@groups for any prime~$p$.

In contrast to these last examples, nilpotency classes are preserved 
when we deal with Aut-bases instead of Inn-bases. Let us denote 
by $\cautg$ the group of all automorphisms of a group~$G$ acting 
trivially on both $Z(G)$ and $G/Z(G)$.  We stick to our terminological 
convention and call any $\cautg$\@basis a CAut-basis of~$G$.   
Now $\cautg=\{1+\eps\st\eps\in\CEndg\}$, where $\CEndg$ is the set of 
all endomorphisms $\eps$ of~$G$ such that $\im\eps\le Z(G)\le 
\ker\eps$.  Thus the property of being a CAut-basis of~$G$ is the 
same as being a CEnd-basis. One can reword
Lemma\rif{AutBasiNormali} as follows: a normal subgroup containing 
the centre is an Aut-basis if and only if it is simultaneously an 
Inn- and a CAut-basis.

\Th\detto{Loryno}  Let the subgroup~$H$ be a CAut-basis of the 
nilpotent group~$G$. Then $Z_n(H)=Z_n(G)\cap H$
for any $n\in\No$. In particular, $H$ and~$G$ have the same  
nilpotency class.

\Pf 
Let $c$ be the nilpotency class of~$G$.  It will be enough to prove 
$Z_n(H)\le Z_n(G)$ for every non-negative integer~$n\le c$.  We will 
argue by induction on~$c-n$, our claim being obvious for $n=c$.  
Thus, suppose $0\le n<c$ and $Z_{n+1}(H)\le Z_{n+1}(G)$.  For every 
integer $i$ such that $0\le i\le n$ let $X_i=[Z_n(H), {}_i G, 
{}_{n-i}H]$.  Our aim is to prove $X_n=1$.  Assume false.  Then we can 
pick the least integer~$i$ between $0$~and~$n$ such that $X_i\not=1$ 
(of course $i>0$) and elements $g_1,g_2,\dots,g_{i-1}\in G$, 
$h_{i+1},h_{i+2},\dots,h_n\in H$ and $a\in Z_n(H)$ such that 
$[a,g_1,\dots,g_{i-1}, G,h_{i+1},\dots,h_n]\not=1$.  The mapping 
$$
\eps:x\in G\mapsto[a,g_1,g_2,\dots,g_{i-1}, x,  
h_{i+1},h_{i+2},\dots,h_n]\in G
$$
belongs to~$\CEndg$, since $a\in Z_{n+1}(H)\le Z_{n+1}(G)$.  Moreover
$H\le\ker\eps$, because
$$
H^\eps=[a,g_1,\dots,g_{i-1}, H, h_{i+1},\dots,h_n]\le X_{i-1}=1.  $$
Since $H$ is a CEnd-basis of~$G$ it
follows~$\eps=0$, that is $[a,g_1,\dots,g_{i-1}, G, 
h_{i+1},\dots,h_n]=1$, contradicting our choice above.  \qed

Neither the analogue of Lemma\rif{abel}$\,(i)$ nor that of 
Theorem\rif{Loryno} holds for hypercentral groups: a term of the upper 
central series of a hypercentral group may be an Aut-basis even if 
its hypercentral length is less than that of the group.
The examples we exhibit to support this statement also have the
feature that the Aut-basis is periodic but the group is not.
For the sake of later reference the construction is slightly 
more general than what is required here.

\Example\detto{noLoryno}
Let $D$ be a locally nilpotent periodic group with infinitely many
non-trivial primary components, none of which is abelian.  Further assume that
$Z(D)$ is reduced.  Let $\pi=\pi(D)$, the set of all primes dividing 
the order of some element of~$D$,  and let $C=\Cr_{p\in\pi}D_p$, the
cartesian product of the primary components of~$D$.  For any $p\in\pi$
let $a_p$ be a noncentral element of~$D_p$.  Consider the subgroup
$A=\Cr_{p\in\pi}\gen{a_p}$ of~$C$ and its element $a=(a_p)_{p\in\pi}$. 
Then $D=\tor C$ and $AD/D$ is a divisible torsion-free abelian group. 
Let $T/D\gen a=\tor(AD/D\gen a)$.  Clearly $T/D\iso\Q$ and $T$ is
locally nilpotent.  Now let $G/D$ be any non-trivial $\pi$\@divisible
subgroup of~$T/D$ (a possible choice is $G=T$).  Then
$C_G(D)=Z(D)=Z(G)$, otherwise we would have $[a^n, D]=1$ for some
$n\in\N$, while $[a^n, D_p]\not=1$ for any $p$~in~$\pi$ not
dividing~$n$.  As $\Hom\(G/D, Z(D)\)=0$ we are now in position to
apply Lemma\rif{AutBasiNormali} and conclude that~$D$ is an Aut-basis
of~$G$.

Now, for each~$p\in \pi$ assume that $D_p$ is nilpotent of class 
greater than~$n_p$ and $a_p\notin Z_{n_p}(D_p)$,
where $\{n_p\st p\in\pi\}$ is unbounded.  
Then $D=Z_\omega(D)=Z_\omega(G)$ and $G=Z_{\omega
+1}(G)$.\qed

Next we discuss subnormal nilpotent subgroups as Aut-bases.

\Lemma\detto{iCentri}
Let $H$ be a normal subgroup of the group~$K$.  Assume  that $H$ is 
contained in the hypercentre~$\bbar Z(K)$ of~$K$. Then  one of the 
following holds: either
\\a) $\Hom\(K/H,Z(K)\cap H\)\not=0$;\qquad or \\b) $Z_i(H)\le Z_i(K)$ 
for all $i\in\No$.

\Pf
For all $j\in\No$ set  $C_j=Z_j(H)$.
Assume that $(b)$ is false. Then the factor $C_{i+1}/C_i$ is not 
central  in~$K$ for some $i\in\No$. Since $H\le \bbar Z(K)$ then 
$C_i<A<B\le C_{i+1}$, where 
$A/C_i=Z(K/C_i)\cap C_{i+1}/C_i$ and $B/A=Z(K/A)\cap C_{i+1}/A$.  Let 
$b\in B\setminus A$.  The mapping $\eps:k\mapsto [k,b]C_i$ from~$K$ 
to~$A/C_i$ is a non-zero homomorphism whose kernel contains~$H$.
For each $j\in\No$ consider the property

\medskip	
\centerline{$\Tondo P_j:$\qquad
\hbox{\bigg\{$\vcenter{\advance \hsize by -5 cm  \line{there exists a 
non-zero homomorphism from $K/H$\hss}  \line{to $C_{j+1}/C_j$ whose 
image is a central factor  of~$K$.\hss}}$}}
\medskip\noindent	
We have just shown that $\Tondo P_i$ holds.  Thus we can pick the  
minimal $j\in\No$ such that $\Tondo P_j$ is true. Suppose~$j>0$.
As~$\Tondo P_{j-1}$ is false the factor $C_j/C_{j-1}$ is central  in~$K$.
Let $\eta$ be a non-zero homomorphism from $K/H$ to $C_{j+1}/C_j$ 
such that $X/C_j:=\im\eta$ is central in~$K$.  There exists $g\in H$ 
such that $[X,g]\not \le C_{j-1}$.  The mapping $$
xC_j\in X/C_j\mapsto [x,g]C_{j-1}\in C_j/C_{j-1} $$
is a non-zero homomorphism whose image is central in~$K$.  Since 
$X/C_j$ is an epimorphic image of~$K/H$ this yields~$\Tondo P_{j-1}$, 
against our choice of~$j$.  Thus $j=0$ and $\Tondo P_0$ holds.  But 
$\Tondo P_0$ is precisely condition~$(a)$ in our statement, so the 
proof is complete.\qed

For the sake of simplifying the next lemma  we introduce a 
piece of notation. 
If $H$ is a subgroup of the group~$K$, we
say that $K$ has property~$(\chi)$ with respect to~$H$ if and only if
$\Hom\(K/H^K,Z(H)\)=0$.

\Lemma\detto{verso il basso}
Let $K$ be a group, and let $H\le \bbar Z(K)$. If 
$K$ has property~$(\chi)$ with respect to~$H$ then also $H^K$ 
satisfies~$(\chi)$ with respect to~$H$.

\Pf
Assume that $N:=H^K$ does not satisfy~$(\chi)$, that is:  
$\Hom\(N/H^N,Z(H)\)\not=0$. Then there exists a proper normal  
subgroup~$S$ of~$N$ such that $H\le S$ and $N/S$ embeds in~$Z(H)$. 
Let $V=S_K$ and  $C/V=(N/V)\cap Z(K/V)$. Obviously $C\cap S=V$, thus 
$C/V$ embeds in~$N/S$ and so in~$Z(H)$.  Now $S\not\n K$ as $N=H^K$, 
hence $C<N$.  As $N\le\bbar Z(K)$ there 
exists $d\in N\setminus C$ such that $dC$ is central in~$K/C$. Let  
$D/V$ be the centralizer of~$dV$ in~$K/V$. Then $N\le D$ (indeed,  
$N/V$ is abelian as $N/S$ is abelian) and $K/D$ is isomorphic to 
$[d,K]V/V$, a non-trivial subgroup of~$C/V$. Since the latter embeds 
in~$Z(H)$, this gives rise to a non-zero homomorphism from $K/N$ 
to~$Z(H)$, in  contradiction to the assumption that $K$
satisfies~$(\chi)$.\qed

\Th\detto{Loryno+}
Let the group $G$ have a nilpotent subnormal subgroup~$H$ contained 
in the  hypercentre and which is an Aut-basis. Then $G$ is nilpotent.

\Pf
Argue by induction on the subnormal defect~$d$ of~$H$ in~$G$. We may 
of course assume~$d>0$.
By Lemma\rif{condHom} we know that $G$ has property~$(\chi)$  with 
respect to~$H$. Lemma\rif{verso il basso} shows that all terms  of 
the standard normal closure series of~$H$ in~$G$  satisfy~$(\chi)$. 
Let $K$ be the term of this series whose defect  is~$d-1$. Then $H\n 
K$ and Lemma\rif{iCentri} gives  $H\le Z_c(K)$ where $c$~is 
the nilpotency class of~$H$. Now  Lemma\rif{abel}$\,(i)$ proves that 
$K$ is nilpotent; by the induction hypothesis also $G$ is nilpotent.\qed

Thus a hypercentral group having a nilpotent subnormal subgroup as 
an Aut-basis is itself nilpotent.  The same does not necessarily 
happen with more general locally nilpotent groups, even if the subgroup 
is normal.  Also the hypothesis that $H$ is an Aut-basis in 
Theorem\rif{Loryno+} cannot be replaced by $H$ being just an Inn-basis.  
This is shown in the next two examples, that can also be compared 
with Proposition\ref{metab}.
We leave open the question whether a hypercentral group with a  
nilpotent subgroup as an Aut-basis must necessarily be  nilpotent. 
% \looseness -1

\Example\detto{27}
For every prime number~$p$ there exists a centre-by-metabelian Fitting
$p$\@group which is not nilpotent and has a normal nilpotent  
subgroup of  class~$2$ (extraspecial, as a matter of fact) as an
Aut-basis.

Indeed, let $P$ be a non-abelian group of order~$p^3$.  Let $A\iso
\pru$, a Pr\"ufer $p$\@group.  Build the standard wreath product
$W=P\wr A$ and let $B=\Dr_{a\in A}P^a$ be its base group.  For every
finite subgroup~$A_0$ of~$A$ the subgroup~$BA_0$ is nilpotent (of the
same class as~$P\wr A_0$), hence $W$ is a Fitting group. 
Set~$Z(P)=\gen c$ and let $Z=Z(B)$.  The elements of~$Z$ are of the
form $c^{n_1 a_1}c^{n_2 a_2}\cdots c^{n_t a_t}$, where $t\in\No$ and
$n_i\in \Z$ and $a_i\in A$ for each~$i$.  Then $K:=[Z,A]$ is the
kernel of the epimorphism~$Z\epi \Z_p$ which maps any such element to
$\sum_{i=1}^t n_i \in\Z_p$.  Let $G=W/K$ and $N=B/K$.  Since
$C_A(B/Z)=1$ it is clear that $C_G(N)=Z(N)$.  Moreover
$Z(N)=Z(B/K)=Z/K$.  Indeed, let $b\in B\setminus Z$.  Then there
exists $a\in A$ such that the image~$b_a$ of~$b$ under the natural
projection $B\epi P^a$ is not central in~$P^a$.  Let $g\in
P^a\setminus C_{P^a}(b_a)$.  Then $[b,g]=[b_a,g]\notin K$, as $P^a\cap
K=1$.  Hence $bK\notin Z(N)$.  Therefore $C_G(N)=Z/K=Z(G)$.  As
$\Hom(G/N, Z(G)\,)\iso \Hom(\C_{p^\infty}, \C_p)=0$,
Lemma\rif{AutBasiNormali} shows that $N$ is an Aut-basis of~$G$. 
That $N$ is extraspecial and $G$ satisfies the requirements of
our claim is now clear.\qed

\Example\detto{22} For every prime number~$p$ there exists a
centre-by-metabelian hypercentral $p$\@group of length~$\omega+1$ with
a class\@$2$ nilpotent subgroup of index~$p$ as an Inn-basis.

Consider first the case $p=2$.
For every $n\in\N$ let $F_n$ be a free group on two 
generators $x_n$~and $y_n$ in the variety of class\@$2$ nilpotent 
groups of  exponent at most~$2^{n+1}$. Thus $F_n=(\gen {x_n}\times  
\gen{c_n})\semid\gen{y_n}$ where $x_n$ and $y_n$ have order~$2^{n+1}$
and $c_n=[x_n,y_n]$ 
has order~$2^n$. Also, $Z(F_n)=\gen{c_n, x_n^{2^n},y_n^{2^n}}$.
Now let $N=\Dr_{n\in\N}F_n$. 
Clearly $N$ has an automorphism~\a of order~2 defined by 
$x_n^\a=x_n^{-1}$  and $y_n^\a=y_n^{-1}$ for all $n\in\N$. Let 
$G=N\semid\gen\a$. Of  course, $G$ is a hypercentral $2$\@group of 
length~$\omega+1$. The equality $C_G(N)=Z(N)$ is 
easily checked.  Now $\a$ centralizes~$Z(N)$, thus $N$ is an 
Inn-basis of~$G$. 

The construction for odd primes~$p$ is rather similar.
% , and we only 
% sketch it. 
For all $n\in\N$  let $F_n$ be
a free group on generators $x_1$, $x_2,\ldots, x_{p-1}$ in the variety
of all nilpotent groups of class at most~$2$ and
exponent dividing~$p^n$. The automorphism~$\alpha_n$ of~$F_n$
defined by $ x_1\mapsto x_2\mapsto \cdots \mapsto x_{p-1}\mapsto 
x_1^{-1}x_2^{-1}\cdots x_{p-1}^{-1}$ has order~$p$. Let
$P_n=F_n\big/[F_n',\alpha_n]$. Now let $N=\Dr_{n\in\N}P_n$ and
$G=N\semid\gen\a$, where $\alpha$ has order~$p$ and acts on each 
direct factor~$P_i$ like~$\alpha_i$. Then $N$ is an Inn-basis of~$G$,
and $G$ is hypercentral, but not nilpotent.
\qed

We close this section with two examples in the class of finite groups. The 
Sylow  $2$\@subgroup of ${\rm SL}(2,3)$ (normal and isomorphic to the 
quaternion group of  order~$8$) is an Aut-basis, obviously hereditary 
since it is a maximal subgroup. More interesting is that the 
Sylow $2$\@subgroups of ${\rm SL}(4,3)$ are End-bases, which proves 
that (even  among finite groups) the existence of a nilpotent 
End-basis does not  necessarily imply nilpotency. 

%%%%%%%%%
%       %
%   %   %
%   %   %
%   %   %
%       %
%%%%%%%%%
\Sezione Hereditary Aut-bases in locally nilpotent groups

The examples above show that a locally nilpotent group~$G$ with a 
nilpotent subgroup~$H$ as an Aut-basis does not need to be 
nilpotent.  We prove in this section that this is the case if we 
strengthen the condition and require that $H$ is a hereditary 
Aut-basis: in this  case $G$ is bound to be nilpotent (of the same 
class as~$H$, by Theorem~\rif{Loryno}). The hypothesis that~$G$ be
locally nilpotent cannot be 
easily  relaxed, as is shown by the polycyclic group in 
Example\rif{ex1}.

Let us start with a couple of simple lemmas. If $G$ is a group and 
$H$  is a subgroup of~$G$, we call a {\it section between $H$ 
and~$G$} a  section $U/V$ of~$G$ such that $H\le V$. If $H$ is a  
hereditary Aut-basis of~$G$ and $X$ is a section between $H$ and~$G$ 
then $\Hom\(X,Z(H)\)=0$ by Lemma\rif{condHom}. 
If $Z(H)\ne 1$ (resp.~if $Z(H)$ has elements of prime order~$p$)
then this shows that~$X$ cannot be infinite cyclic (resp.~of order~$p$).
Thus:

\Lemma\detto{no-sect} Let the subgroup $H$ be a hereditary Aut-basis 
of the group~$G$. If $Z(H)\not=1$ then:
\\i) there is no infinite cyclic section between $H$ and $G$;
\\ii) for every prime 
$p\in\pi\(Z(H)\)$ there is no section of order~$p$ between $H$ and $G$.

\Lemma\detto{semid}Let $G=N\semid K$ be a semidirect product. If $K$ 
is  a hypercentral Aut-basis of~$G$, then $G=N\times K$ and $|N|\le 
2$.

\Pf
Let $C=C_K(N)$. Argue by means of contradiction and assume $C<K$.  
Then $K/C$ has a nontrivial central element~$xC$. The conjugation 
automorphism induced by $x$ on~$N$ commutes with those induced by  
elements of~$K$, hence the mapping defined by $ak\mapsto a^xk$ (for  
$a\in N$ and $k\in K$) is a non-trivial automorphism of~$G$ which  
acts trivially on~$K$. This is a contradiction, since $K$ is an Aut-basis 
of~$G$. Therefore $C=K$, that is: $G=N\times K$. Again the fact that $K$ 
is an  Aut-basis yields $\aut N=1$, hence $|N|\le 2$.\qed

The next lemma and the corollary following it are doubtless well 
known; we include proofs since we could not find suitable references 
to them.

Let $H$ be a subgroup of the locally nilpotent group~$G$ and let $\pi$ be 
a set of primes. We write~$I_{G,\pi}(H)$ for the {\it 
$\pi$\@isolator} of~$H$ in~$G$, i.e., the subgroup consisting of all 
elements~$x$ of~$G$~such that $x^n\in H$ for some $\pi$\@number~$n$. 
According to standard  terminology  $H$ is {\it $\pi$\@isolated} 
in~$G$ if and only if $H=I_{G,\pi}(H)$. Also, $\pi'$ denotes the set 
of all primes not in~$\pi$. 
If $\pi$ is the set of all 
primes we write $I_G(H)$ for~$I_{G,\pi}(H)$, the isolator of $H$ 
in~$G$.   We shall make use of the following two well-known facts:

\smallskip\noindent
(1) if $H$ is hypercentral of length~$\a$ and $G=I_{G,\pi}(H)$ has
trivial $\pi$\@component, then also $G$ is hypercentral of length~$\a$
(indeed, by adapting the proof of~\[EN], Lemma~4.8, for instance, one
shows that $Z_\beta(G)$ is $\pi$\@isolated in~$G$ and contains
$Z_\beta(H)$ for every ordinal~$\beta$.  Hence $Z:=Z_\a(G)\ge H$ and
so $Z=I_{G,\pi}(Z)\ge I_{G,\pi}(H)=G$);

\smallskip\noindent
(2) if $G$ is a finitely generated nilpotent group and $H\le G$, then 
$|G:H|$ is finite if and only if $G/HG'$ is finite.

\Lemma\detto{max}  Let $x$ be an element of infinite order in the 
locally  nilpotent group~$G$. Let $H$ be maximal among the subgroups 
of~$G$  subject to the condition $H\cap \gen x=1$. Then $x$ 
normalizes~$H$.

\Dim
Let \Sy be a local system of finitely generated  subgroups of~$H$. 
For every $K\in\Sy$ let  $K^*=K\gen{K,x}'$. Since 
$\big|\gen{K,x}:K\big|$ is infinite the  remark immediately 
preceding this lemma shows that  the cyclic factor $\gen{K,x}/K^*$ 
is infinite. 	 Let $\Sy^*=\{K^*\st K \in\Sy\}$. The mapping 
$K\in\Sy\mapsto  K^*\in\Sy^*$ preserves inclusion, hence $\Sy^*$ is 
a local system of  subgroups of $H^*:=\gen {K^*\st K \in\Sy}$. Of 
course $x$  normalizes~$H^*$ and $H\le H^*$. Suppose $H<H^*$. Then, 
by maximality  of~$H$, there exists $n\in\N$ such that $x^n\in H^*$. 
But then $x^n\in  K^*$ for some $K^*\in \Sy^*$. This is a 
contradiction, since  $\gen{K,x}/K^*$ is infinite cyclic. Therefore  
$H=H^*$ is normalized by~$x$.\qed

\Cor\detto{dicotomia}
Let $H$ be a subgroup of the locally nilpotent group~$G$. Then  
either $G=I_G(H)$ or there is an infinite  cyclic section between 
$H$~and~$G$.

\Dim
If $G$ is not the isolator of~$H$ in~$G$ there exists in~$G$ an 
element $x$  of infinite order such that $H\cap\gen x=1$. By Zorn's 
Lemma there is  a subgroup $M$ of $G$ maximal with respect to $H\le 
M$ and $M\cap\gen  x=1$. It follows from Lemma~\ref{max} that 
$\gen{M,x}/M$ is an  infinite cyclic section between $H$ and~$G$.\qed

\Lemma\detto{ipercbasi}
Let $H$ be a hypercentral subgroup of 
the locally nilpotent group~$G$. Let $\pi=\pi(H)$. If $H$ is a  
hereditary Aut-basis of~$G$, then:
\\i)if $U/V$ is a periodic section of~$G$ such that $\pi(U/V)\subseteq\pi$, then $\gen{H,U}=\gen{H,V}$;
\\ii) $G=I_{G,\pi'}(H)$;
\\iii)$\tor H$ is the $\pi$\@component of~$\tor G$; the 
$\pi'$\@component of~$\tor G$ has order at most~$2$.

\Pf
$(i)$\enskip Assume
$X:=\gen{H,V}<\gen{H,U}$ and let $x\in U\setminus X$.  There is a 
maximal subgroup~$M$ of~$Y:=\gen{X,x}$ containing~$X$.  Since $G$ is 
locally nilpotent and $x^n\in V\le X$ for some $\pi$\@number~$n$, the 
order of~$Y/M$ is a prime belonging to~$\pi$. But  $Y/M$ is a section 
between $H$ and~$G$. This and Lemma~\ref{no-sect}\thinspace$(ii)$ 
yield a contradiction. 
\\ii) Equality $G=I_G(H)$ immediately follows from  
Lemma~\ref{no-sect}\thinspace$(i)$~and Corollary~\ref{dicotomia}. On 
the  other hand, what we have just proved as part~$(i)$ ensures that 
$H$ is $\pi$\@isolated in~$G$. It is then clear that 
$G=I_{G,\pi'}(H)$ holds.

\\iii) Let $P$ be the $\pi$\@component of~$\tor G$. Then 
$\gen{H,P}=H$ by part~$(i)$, hence $P\le H$.  This shows that $\tor 
H$ is the $\pi$\@component of~$\tor G$. Finally,  if $N$ is the 
$\pi'$\@component of~$\tor G$, then $H$ is an Aut-basis of $\gen 
{H,N}=N\semid H$, so $|N|\le 2$  by Lemma~\ref{semid}.\qed

Now we are in position to prove the main result of this section.

\Th\detto{Teorema Pasqualino} 
Let the locally nilpotent group~$G$ have a 
nontrivial hypercentral subgroup~$H$ 
as a  hereditary Aut-basis.  Then $G$ is hypercentral of  the same 
hypercentral length as~$H$. In particular, if $H$ is nilpotent  then 
also $G$ is nilpotent.

\Pf
Let $H$ be hypercentral of length~$\a$. At the expense of 
replacing~$H$ by $(\tor  G)H$---which is still hypercentral of 
length~$\a$ by  Lemma\rif{ipercbasi}\thinspace$(iii)$---we may assume 
$\tor G\le  H$. Set $\pi=\pi(H)$. Then $G$ has trivial 
$\pi'$\@component and  $I_{G,\pi'}(H)=G$ by 
Lemma\rif{ipercbasi}\thinspace$(ii)$. Thus $G$ is  hypercentral of
length~$\a$, as required.\qed

We recall from Example\rif{noLoryno} that even if $G$ is hypercentral 
but the  Aut-basis~$H$ is not hereditary it may happen that the 
hypercentral  lengths of~$H$ and~$G$ differ.

\ignora{ %%%%%%%%%%%%%%%%%%%%%%%
As a final remark we note a special case of  Lemma~\ref{ipercbasi} 
$(iii)$. The latter shows that, apart from trivial instances, no 
hypercentral subgroup may be a
hereditary Aut-basis in a periodic locally nilpotent group.

\Cor\detto{ipercentr. periodici ereditarie}
Let $G$ be a periodic locally nilpotent group. Then a  
hypercentral subgroup~$H$of~$G$ is a hereditary Aut-basis if and only 
if either $G=H$ or $2\notin \pi(H)$ and $G=H\times K$ with $|K|=2$.

\Dim
The `only if' part is in Lemma~\ref{ipercbasi}; the `if' part is
obvious as $K$ is characteristic in~$G$ and $\aut K=1$, if $H\not=
G$.\qed
%%%% IGNORATO %%%%%%%%%%%
}

As a final remark we note that Lemma~\ref{ipercbasi} shows that, apart
from trivial instances, no hypercentral subgroup may be a hereditary
Aut-basis in a periodic locally nilpotent group.  It is obvious (see
Lemma\rif{no-sect}\thinspace$(i)\,$) that a locally nilpotent group is
periodic if it has a periodic subgroup with non-trivial centre which
is a hereditary Aut-basis.  So we have:

\Cor\detto{ipercentr. periodici ereditarie}
Let $G$ be a locally nilpotent group. Then a  periodic
hypercentral subgroup~$H$of~$G$ is a hereditary Aut-basis if and only 
if either $G=H$ or $2\notin \pi(H)$ and $G=H\times K$ with $|K|=2$.

\Dim
The necessity of the condition follows from
Lemma~\ref{ipercbasi}\thinspace$(iii)$ and the last observation.  Its
sufficiency is immediate, as $K$ is characteristic in~$G$ and $\aut K=1$,
if $H\not=G$.\qed

If the hypotheses on~$H$ are relaxed, it is possible to obtain less
obvious examples of hereditary Aut-bases in locally nilpotent groups. 
Indeed, in a torsion-free abelian group every maximal independent subset is a
hereditary Aut-basis.  Also, if $D$ is a locally nilpotent periodic
group with trivial centre and infinitely many nontrivial primary
components, the construction in Example\rif{noLoryno} embeds~$D$ as a
hereditary Aut-basis (and the torsion subgroup) of a non-periodic
locally nilpotent group; this follows by a direct application of 
Lemma\rif{centroidentico}.

\Sezione Finite Aut-bases and outer automorphisms

If $G$ is a periodic countably infinite 
locally nilpotent group then $|\autg|=2^{\aleph_0}$.
This was proved by O.~Puglisi in~\[OP]. According to 
Kueker's theorem quoted in the introduction this result is equivalent 
to the following: a countable periodic locally nilpotent group has a 
finite Aut-basis if and only if it is finite. After rewording 
Puglisi's theorem in this form, inspection of the 
proofs in~\[OP] reveals that the countability hypothesis can be 
disposed of. Here we aim at further extensions of this theorem, 
by relaxing the periodicity hypothesis and also by weakening the 
requirement of being an Aut-basis.

Let $G$ be a group.  With respect to subgroups $\G$ of~$\autg$ it can 
be of interest to consider $C_{\autg}(\G)$\@bases.  This is because of 
the following: assume $G\n E$ for some group~$E$ and $E=GK$ for some 
$K\le E$. Then there exists no non-trivial automorphism of~$E$ which 
acts trivially on~$K$ and normalizes~$G$ if and only if $H:=K\cap G$ 
is a $C_{\autg}(\G)$\@basis, where \G is the group of all 
automorphisms induced on~$G$ by conjugation by elements of~$K$.  
Hence $H$ is a $C_{\autg}(\G)$\@basis if $K$ is an 
Aut-basis of~$E$, and the two properties are equivalent if $G$ is 
characteristic in~$E$. 

\ignora{ %%%%%%%%%%%%%%%%%%%%%%%
% Hence the property that $H$ is a $C_{\autg}(\G)$\@basis is weaker 
% than $H$ being an Aut-basis of~$G$ and also than $K$ being an 
% Aut-basis of~$E$, but it is equivalent to the latter if $G$ is 
% characteristic in~$E$. 

% Also note that in such a situation $H$ is 
% $\G$\@invariant. 
% For this reason we will state  consider  
% 
%%%%% IGNORATO %%%%%%
}%%%%%%%

We will make use of the following variation of Lemma~4 of~\[OP].

\Lemma\detto{Philemma}
Let $G$ be a locally nilpotent group,~$\pi$~a finite set of  primes, $n$ the 
product of the primes in~$\pi$. Assume that a $\pi$\@subgroup~$H$ 
of~$G'G^n$ is an Inn-basis of~$G$. If either  $G$ is
hypercentral or $H$ is finite, then $G$ is abelian.

\Dim
Assume first that $G$ is hypercentral. Suppose that $G$ is not abelian. 
Then $H\not\le Z(G)$, as $H$ is an Inn-basis of~$G$.  Therefore
$Z_2(G)\big/Z(G)$ has non-trivial $p$\@component for some $p\in\pi$. 
Let $xZ(G)$ be an element of order~$p$ in $Z_2(G)\big/Z(G)$.  Then
$G/C_G(x)\iso[G,x]$ is an elementary abelian $p$\@group, so that
$G'G^p\le C_G(x)$.  Since $G'G^p$ is an Inn-basis of~$G$ this yields
$x\in Z(G)$, a contradiction. Thus $G$ is abelian in this case.

Assume now that $H$ is finite.
Suppose again that $G$ is not abelian. 
Then from $H\le G'G^n$ it follows 
that $G$ has a finitely generated non-abelian subgroup~$X$ 
such that $H\le X'X^n$. Since $X$ is nilpotent and $H$ is an Inn-basis
of~$X$ this is impossible by what was proved in the previous case.\qed

\Lemma\detto{centralizzanti}
Let \G be a group of automorphisms of the group~$G$ 
such that the external semidirect product $E=G\semid \G$ is locally 
nilpotent. Let the $\G$\@invariant subgroup $H$ of~$G$ be a  
$C_{\autg}(\G)$\@basis. 
If either \G is finitely generated or $E$ is hypercentral, 
then $C_E(H)=C_E(G)$; in particular $H$ is an 
Inn-basis of~$G$.

\Dim
Let $\Delta$ be the group of those automorphisms of~$G$ induced by 
conjugation by elements of~$C_E(H)$.  Since $H$ is $\G$\@invariant 
$\Delta$ is normalized by~$\G$, and $\Delta\G\iso C_E(H)\G/C_E(G)$ is 
locally nilpotent.  Now $C_\Delta(\G)=1$, because $H$ is a 
$C_{\autg}(\G)$\@basis and $[\Delta,H]=1$. By hypothesis either
\G is finitely generated or $\Delta\G$ is hypercentral. In both cases
from $C_\Delta(\G)=1$ it follows $\Delta=1$, 
which amounts to saying $C_E(H)=C_E(G)$.\qed

The special case dealt with in the next lemma is of great relevance for
the proofs that follow.

\Lemma\detto{abelian case}
Let $\pi$ be a set of primes, and let \G be a periodic hypercentral group of
automorphisms of the $\pi$\@divisible abelian group~$G$ such that the
external semidirect product $E=G\semid \G$ is locally nilpotent. 
Assume that $G$ has a $\pi$\@subgroup~$H$ of finite exponent which is
a $C_{\autg}(\G)$\@basis.  Then both $\G$ and the $\pi$\@component
of~$G$ are trivial, and $|G|\le 2$.

\Dim
Let $e=\exp H$ and~$\psi=\pi(H)$.
At the expense of replacing $H$ by $G[e]=\{g\in G\st g^e=1\}$
we may assume that $H$ is $\G$\@invariant.
The $\psi$\@component~$P$ of~$G$ is
divisible, so $G=P\times Q$ for some subgroup~$Q$.

If $\G=1$ then $H$ is an Aut-basis of~$G$.
For any $\a\in\aut  Q$ we can define 
an automorphism of~$G$ which acts on~$Q$ like~$\a$ and maps every 
element of~$P$ to its $(e+1)$\@th power. This automorphism 
centralizes~$H$, hence it is the identity. Thus $\aut Q=1$ and  
$P^{e}=1$. It follows $|Q|\le 2$ and $P=1$, since $P$  
is divisible; by the same reason the $\pi$\@component of~$G=Q$ is trivial.
Thus all we have to prove is that \G is trivial. 

Assume $\G\not=1$. Let $C=C_\G(H)$. If $C\not=1$ then
$C_C(\G)\not=1$, because \G is hypercentral~and $C\n\G$. 
This is impossible, as $H$ is a $C_{\autg}(\G)$\@basis. Thus $C=1$, so
that \G is isomorphic to a periodic group of
automorphisms of~$H$. As $H\semidirect \G$ is locally nilpotent it follows 
that \G is a
$\psi$\@group.  This also gives $[G,\G]\le P$
(indeed, let $X$ be a finitely generated subgroup of $G/P$, and let 
$\gamma\in\G$.  Then
$X^{\gen\gamma}$ is a finitely generated abelian group with trivial
$\psi$\@component on which the $\psi$\@group~$\gen\gamma$
acts nilpotently, hence $[X^{\gen\gamma},\gamma]=1$). 
With respect to the decomposition $G=P\times Q$ any
automorphism~$\gamma$ of~$G$ acting trivially on~$G/P$
(hence any element of~$\G$) can be
represented by the matrix~$\smatrix(\a,0,\eps,1)$, where
\a is the automorphism induced by~$\gamma$ on~$P$ and 
$\eps\in \Hom(Q,P)$ is defined by $x^\gamma=x^\eps x$ for 
every $x\in Q$. 
If $\delta$ is another such automorphism of~$G$, represented
by~$\smatrix(\beta,0,\eta,1)$, then $\gamma$~and~$\delta$
commute if and only if the following conditions hold:
$$
\a\beta=\beta\a\qquad{\rm and}\qquad \eps(\beta-1)=\eta(\a-1).
\eqno(*)
$$
As $\G\not=1$ there exists a nontrivial element $\zeta\in Z(\G)$. 
Assume that $\zeta$ is represented by~$\smatrix(\a,0,\eps,1)$.  Let
$\theta$ be the automorphism of~$G$ represented
by~$\smatrix(\sigma,0,\tau,1)$, where $\sigma=1+e(\a-1)$ and
$\tau=e\eps$.  To make sure that this definition is consistent we have
to check that such~$\sigma$ is indeed an automorphism.  It is clear
that $\sigma$ is a monomorphism, since it acts trivially on the socle
of~$P$.  We have to prove that~$\sigma$ is epic.  Every $x\in P$
belongs to a finite $\a$\@invariant subgroup~$X$ of~$P$, as \a is
periodic.  Certainly $X^\sigma\le X$, hence $X^\sigma= X$
because~$\sigma$ is a monomorphism, so $x\in P^\sigma$.  Therefore
$\sigma\in\aut P$.  It is obvious that~$\theta$ acts trivially on~$H$;
we claim that~$\theta$ centralizes~$\G$.  To prove this claim, let
$\gamma\in\G$ and assume that~$\gamma$ is represented
by~$\smatrix(\beta,0,\eta,1)$.  Since $[\gamma,\zeta]=1$ equalities as
in~$(*)$ hold.  These yield $\sigma\beta=\beta\sigma$ and
$$
\tau(\beta-1)=e\eps(\beta-1)=e\eta(\a-1)=\eta e(\a-1)=\eta(\sigma-1).
$$
Therefore $[\theta,\gamma]=1$ and our claim is established.
It follows $\theta=1$, as $H$ is a $C_{\autg}(\G)$\@basis. Hence
$\sigma=1$ and $\tau=0$, which in turn give 
$\a=1$ and $\eps=0$, because $P$ and~$Q$ are $\psi$\@divisible.
This means $\zeta=1$, against our choice. This contradiction 
shows $\G=1$, so the proof is complete.\qed

\Th\detto{finito}
Let \G be a finitely generated group of automorphisms of the group~$G$ 
such that the external semidirect product $E=G\semid \G$ is locally 
nilpotent. Assume that $G$ has a finite subgroup~$H$ which is a 
$C_{\autg}(\G)$\@basis. Then $E$ is finite.

\Dim
It is clear that the
subgroup~$H^\G$ is still finite (and a $C_{\autg}(\G)$\@basis). 
Hence---by substituting $H^\G$ for~$H$ if necessary---we may assume
that $H$ is $\G$\@invariant.  Then \G is finite since $H$ is finite
and $C_\G(H)=C_\G(G)=1$ by Lemma~\ref{centralizzanti}.

Suppose that the theorem is false and choose a counterexample in 
which $H$ is $\G$\@invariant of the minimal possible order.
If $H=1$ then $\G=1$ and $|G|\le 2$, hence $H\not=1$.

Assume that $G$ has no maximal subgroups whose index belongs
to~$\pi:=\pi(H)$.  Then $G=G'G^n$, where $n$ is the product of all
primes in~$\pi$, hence $G$ is abelian by Lemma\rif{Philemma} and
Lemma\rif{centralizzanti}. Moreover $G$ is also $\pi$\@divisible, thus
Lemma\rif{abelian case} shows that $E$ is finite.

Therefore $G$ has some subgroup of index~$p$, for some
$p\in\pi$.  Since \G is finite $G$ must then have some $E$\@invariant
proper subgroup of finite $p$\@power index, so it has a maximal
$E$\@invariant proper subgroup, say~$M$.  As $E$ is locally nilpotent
$[G,\G]\le M$ and $\big|G/M\big|=p$.  Suppose first $H\le M$.  Since
$H\semid\G$ is nilpotent the centralizer~$A$ of~$\G$ in the
$p$\@component of~$Z(H)$ is nontrivial, hence there exists a non-zero
homomorphism $\eps:G\rightarrow A$ with kernel~$M$.  Also, $A\le Z(G)$
by Lemma~\ref{centralizzanti}.  Then $g\mapsto gg^\eps$ defines a
non-trivial automorphism of~$G$ centralizing~$\G$ and~$H$.  This is
impossible, as $H$ is a $C_{\autg}(\G)$\@basis.  Hence $H\not\le M$,
so $G=MH$.  Let $\bbar\G$ and $\Delta$ be the groups of the
automorphisms of~$M$ induced by~\G and (by conjugation) by~$H$
respectively.  Clearly $\bbar\G$ normalizes~$\Delta$, the product
$\bbar\G\Delta$ is finite and the external semidirect product
$M\semid\bbar\G\Delta$ is locally nilpotent.  Every automorphism of~$M$
centralizing~$\bbar\G\Delta$ and acting trivially on $H\cap M$ extends
to an automorphism of~$G$ centralizing~\G and acting trivially on~$H$. 
Since $H$ is a $C_{\autg}(\G)$\@basis the latter automorphism is bound to be the
identity.  This proves that $H\cap M$ is a $C_{\aut
M}(\bbar\G\Delta)$\@basis.  Also, it is $\bbar\G\Delta$\@invariant.  By
minimality of~$|H|$ we conclude that~$M$ is finite, hence the same
holds for~$E$, as wished.\qed

The statement of Theorem\rif{finito} can be further improved by
showing that the subgroup $H^\G$ is close to being equal to~$G$.  This
is not needed for the next two corollaries and will be proved as
a special case of Theorem\rif{expfinito}.

\Cor\detto{per-by-fg}
A periodic-by-(finitely generated) locally nilpotent group has a 
finite Aut-basis if and only if it is finitely generated.

\Pf
Let $G$ be a locally nilpotent group, and let $T=\tor G$.  Assume that 
$G/T$ is finitely generated and $G$ has a finite Aut-basis~$X$.  
There exists a finite subset~$Y$ of~$G$ such that $G=\gen {T,Y}$.  
Let $H=\gen {X,Y}$.  Then $H$ is a finitely generated Aut-basis 
of~$G=TH$.  Because of the argument sketched at the beginning of this 
section this means that $T\cap H$ is a $C_{\aut T}(\G)$\@basis, where 
\G is the group of the automorphisms induced by conjugation by~$H$ 
on~$T$.  Thus Theorem\rif{finito} proves that $T$ is finite, hence 
$G$ is finitely generated.\qed

\Cor\detto{outer}
Let $G$ be a countably infinite locally nilpotent group. If\/ $\tor G$  
is infinite and $G/\tor G$ 
is finitely generated then $|\autg|=2^{\aleph_0}$.

\Pf
Suppose $|\autg|\not=2^{\aleph_0}$. Then $G$ has a  finite Aut-basis 
by Kueker's theorem. By Corollary\rif{per-by-fg} then $G$ is finitely 
generated, i.e., $\tor G$ is finite.\qed 

Thus periodic-by-finitely generated locally nilpotent groups which 
are countable but not finitely generated have (uncountably many) outer 
automorphisms.  It appears to be not known whether all infinite 
finitely generated nilpotent groups have outer automorphisms---of 
course the automorphism group of any such group is countable anyway.  
A detailed discussion on this and related topics is in~\[D].
How much the countability hypotesis is decisive for the existence of 
outer automorphisms, even in the periodic case, has been shown 
% in~\[ST] and~\[DG],
by S.~Thomas in \[ST] and by Dugas and G\"obel in~\[DG]. Indeed, the latter
prove that for every infinite cardinal~$\lambda$ such that 
$\lambda=\lambda^{\aleph_0}$ and for every prime~$p$
there exists a complete locally nilpotent
$p$\@group of cardinality~$\lambda^+$.

A better version of Theorem\rif{finito}, more consistent with the
general setting of this paper, would be obtained if the following
question could be answered in the positive: is it possible to replace
the hypothesis that $H$ is finite with the hypothesis that $H$ is
finitely generated still getting that $E$ is nilpotent?  (Of course,
as already remarked even for the rational group, one cannot hope to
prove that~$E$ is finitely generated.)  We have been unable to settle
this conjecture.  If true this would give as a special case, for
$\G=1$, that every locally nilpotent group with a finite Aut-basis is
nilpotent.  Again by Kueker's theorem, this would further imply that
every non-nilpotent countable locally nilpotent group would
have~$2^{\aleph_0}$ (outer) automorphisms.  The well-known example by
A.E.~Zaleski\u\i~\[Z] of a countably infinite nilpotent group with no
outer automorphisms shows that the analogous statement fails to hold
for nilpotent groups.  On the opposite side, as remarked by S.~Thomas
(\[ST], Proposition~3), non-trivial centreless locally nilpotent
groups clearly have no finite Aut-basis (actually, no finite
Inn-basis), hence if such a group~$G$ is countable then
$|\autg|=2^{\aleph_0}$.

\smallskip
In connection with Theorem\rif{finito}, 
it is worth mentioning that a non-periodic nilpotent (or even abelian) group may well 
have a periodic subgroup as an Aut-basis. Probably the most obvious 
example is a cartesian product~$C$ of infinitely many  groups of pairwise 
different prime orders, whose torsion subgroup~$D$ (i.e., the 
corresponding direct product) is easily seen to
be an End-basis, as $\Hom(C/D,C)=0$. Other examples are
those constructed in Example\rif{noLoryno}: in the notation used 
there, if $D$ is
nilpotent then also $G$ is nilpotent. 

One can ask whether it is possible to obtain similar examples where 
the Aut-basis has finite exponent. At least for nilpotent groups, a proof
similar to that of Theorem\rif{finito} shows that the 
answer is negative. To state this fact in a more general setting
we point at the following two properties relative to
groups~$G$:
\def\prM{(\Tondo M)}
\def\prN{(\Tondo N)}
$$
\leqalignno{
&\hbox{\rm every proper subgroup of~$G$ is contained in a
maximal subgroup of~$G$;}&\prM:\cr
&\hbox{\rm if $H$ is a proper subgroup of~$G$ then $HG'<G$.}&\prN:\cr}
$$

As is well known, property~$\prN$ is a form of generalized nilpotency, in the
sense that nilpotency implies~$\prN$ and it is equivalent to~$\prN$ 
for finite groups. Nilpotent groups of finite exponent satisfy~$\prM$.
Indeed, $\prM$~and~$\prN$ are equivalent conditions for a locally
nilpotent group of finite exponent.

\ignora{%%%%%%%%%%%%%%%%%%%%%%%%%%%%%%%
% We stress that $\prM$~and~$\prN$
% can fail to hold for such a group. This is shown, for instance, by the
% locally nilpotent perfect groups of finite exponent constructed
% in~\[VW].
%%%%%%%%%%%%%%%%%%%%%%%%%%%%%%%%%%%%%%%%

\Th\detto{expfinito}
Let $G$ be a hypercentral group all whose quotients of finite
exponent satisfy~$\prN$. 
Let \G be a finite group of automorphisms of~$G$ 
such that the external semidirect product $E=G\semid \G$ is
locally nilpotent. Assume that $G$ has a
subgroup~$H$ of finite exponent~$e$
which is a $C_{\autg}(\G)$\@basis.
Further, assume that either $H$ is $\G$\@invariant or $G$ is nilpotent.
Then $\exp G$ is finite.
Moreover, if $H$ is $\G$\@invariant then either 
$G=H$ or $e$ is odd and $G=H\times N$ with $|N|=2$.

}

\Th\detto{expfinito}
Let $G$ be a hypercentral group all whose quotients of finite
exponent satisfy~$\prN$. 
Let \G be a finite group of automorphisms of~$G$ 
such that the external semidirect product $E=G\semid \G$ is
locally nilpotent. Assume that $G$ has a
subgroup~$H$ of finite exponent~$e$
which is a $C_{\autg}(\G)$\@basis. Then:
\\i) if $G$ is nilpotent of class~$c$ then $\exp G$ is finite and divides
the least common multiple of~$2$ and~$e^c$;  
\\ii) if $H$ is $\G$\@invariant then either 
$G=H$ or $e$ is odd and $G=H\times N$ with $|N|=2$.

\Dim
If $G$ is nilpotent of class~$c$ then 
$H^\G$ has finite exponent dividing~$e^c$. Therefore~$(i)$ is a
consequence of~$(ii)$ and it will be enough to prove the latter.

So, suppose $H^\G=H$.  Let $N=G^{e^2}$.  Assume $NH<G$.  Since $G/N$
satisfies~$\prM$ there is a maximal subgroup of~$G$ containing~$NH$. The
index of this subgroup in~$G$ is a prime~$p$ belonging to~$\pi:=\pi(H)$. As
$NH$ is $\G$\@invariant and \G is finite, it follows that
$NH$---and hence~$H$--- is
contained in an $E$\@invariant subgroup of index~$p$. Since $H\semid 
\G$ is hypercentral (actually, $E$ is hypercentral, as a
hypercentral-by-finite locally nilpotent group) this
can be excluded by the same argument 
as in the proof for Theorem\rif{finito}.  Thus $NH=G$.  Hence $G/N\iso
H/H\cap N$ has exponent
at most~$e$. Therefore $N=G^e$, so $N=N^e$.  Now $H\cap N$ is a
$C_{\aut N}(\bbar\G\Delta)$\@basis, where $\bbar\G$ and $\Delta$ are the
groups of automorphisms of~$N$ induced by~\G and (by conjugation) by~$H$
respectively. Also, $N\semid \bbar\G\Delta$ is hypercentral, thus 
we may apply Lemma\rif{centralizzanti}
and Lemma\rif{Philemma} to obtain that
$N$ is abelian and $\pi$\@divisible. Now Lemma\rif{abelian
case} implies that both the $\pi$\@component of~$N$ and~$\G$ are trivial,
and $|N|\le 2$. The conclusion is now clear. \qed

As for Corollary\rif{ipercentr. periodici ereditarie}, a converse 
statement trivially holds: if $G=H\times N$ as in the last part of 
Theorem\rif{expfinito} then $H$ (the odd component of~$G$) is a
characteristic Aut-basis of~$G$.

Two further remarks about Theorem\rif{expfinito}
are in order.  Firstly, in the hypotheses of the
theorem it is possible that the subgroup~$H$ differs greatly
from~$G$ if it is not $\G$\@invariant.  Indeed, let
$G=\gen x\semid \gen a$ be a dihedral $2$\@group of order at
least~$8$.  Let $\G=\autg$.  Then $G\semid\G$ is a $2$\@group, and
it is easily seen that $\gen a$ is a $Z(\G)$-basis, that is a
$C_{\autg}(\G)$\@basis.

Secondly, if $\G=1$ the $\G$\@invariance condition becomes trivial, so
Theorem\rif{expfinito} gives an explicit description of the Aut-bases 
of~$G$, showing that only the obvious cases arise. Thus we have:

\Cor\detto{nilp exp finito esplicite}
Let $G$ be a periodic locally nilpotent group all of whose primary
components are nilpotent, and let $H\le G$.
Assume that every 
primary component of~$H$ has finite exponent. Then $H$ 
is an Aut-basis of~$G$ if and only if either
$G=H$ or $H$ is the $2'$\@component of~$G$ and
the $2$\@component of~$G$ has order~$2$.

\ignora{\Dim
The statement follows immediately from Theorem\rif{expfinito}, once 
observed that $H$ is  an Aut-basis of~$G$ if and only if the $p$\@component
of~$H$ is an Aut-basis of the $p$\@component of~$G$
for every prime~$p$.\qed}

If $p$ is a prime and $G$ is a direct product of cyclic $p$\@groups 
of unbounded orders then $G$ has a subgroup~$H$ such that $G/H$ is a 
Pr\"ufer group. Now $\Hom(G/H,G)=0$, so $H$ is an End-basis of~$G$.
This shows that
the hypothesis that the primary components of~$H$
have finite exponent is crucial in Corollary\rif{nilp exp finito
esplicite}. Also, the group $G=N\semid A$ 
in Example\rif{27} has the normal subgroup~$N$ 
of finite exponent as an Aut-basis, yet it is of infinite exponent. 
It is also possible to modify this example in order to make it into a 
non-periodic (but still centre-by-metabelian and Fitting) group
with the same Aut-basis. To 
this aim let $A_1$ be a group isomorphic to $\Q\kern1pt_p=\{np^m\mid 
n,m\in\Z\}\le\Q$ and let ~$\eps:A_1\shepi A$.
Let $G_1=N\semid A_1$, where $A_1$ acts on~$N$ 
via~$\eps$ and the action of~$A$. 
%
% Now,  $????$ hence $K:=\ker\eps=A_1\cap Z(G_1)$
% is characteristic in~$G_1$.
% If $\a\in C_{\autg}(N)$ then $[G_1,\alpha]\le K$, because
% $N$ is an Aut-basis in $G\iso G_1/K$. Thus $x\in G_1\mapsto 
% x^{-1}x^\alpha\in K$ is a homomorphism. However $K\iso\Z$, so
% $\Hom(G_1,K)=0$. This shows $\alpha=1$. 
Then $C_{G_1}(N)=Z(N)\ker \eps=Z(G_1)$. Let $\G=C_{\aut G_1}(N)$. 
By the Three Subgroup Lemma $[G_1,\G]\le Z(G_1)$. The only
automorphism of $G_1/N$ that acts trivially on~$G_1/NZ(G_1)$ is the 
identity. Hence $\G$ centralizes ($N$ and) $G_1/N$ and can be embedded in
$\Hom\(G_1/N,Z(N)\)\iso \Hom\(\Q\kern1pt_p,\C_p\)=0$. 
Therefore $N$ is an Aut-basis 
of~$G_1$, as required.
% 
% 
% %let $G_1=N\semid A_1$, where
% replace the Pr\"ufer group~$A$ with $A_1\iso \Q_p=\{np^m\mid 
% n,m\in\Z\}\le\Q$, and 
% 
% 
% a non-periodic
% $p$\@divisible abelian  group~$A_1$ such that $A$ is the image of an 
% epimorphism~$\eps$ defined on~$A_1$, and let $A_1$ act on~$N$ 
% via~$\eps$ and the action of~$A$. By the same argument as in the 
% original example $N$ will still be an Aut-basis.
% this aim let $G_1=N\semid A_1$, where
% replace the Pr\"ufer group~$A$ with a non-periodic
% $p$\@divisible abelian  group~$A_1$ such that $A$ is the image of an 
% epimorphism~$\eps$ defined on~$A_1$, and let $A_1$ act on~$N$ 
% via~$\eps$ and the action of~$A$. By the same argument as in the 
% original example $N$ will still be an Aut-basis.

\bigskip\bigskip
\centerline{\caps Acknowledgments}
\smallskip\noindent
The authors wish to thank the referee for a number of useful comments;
in particular for suggesting the statements of
Lemma\rif{abel}$\;(iii)$ and Proposition\rif{metab}, and for bringing
to their attention the paper~\[DG].
\bigskip
% and for suggesting to
% consider metabelian groups .

\begingroup
\references 
[B] \author R. Baer 
\title Gruppen mit abz\"albaren Automorphismengruppen \journal 
Hamburger Math. Einzelschriften (Neue Folge), \vol{\rm  Heft}~2 \year 
1970 \pag. 1--122.

% %[BW] \author J. Buckley \and J. Wiegold  %\title On the number of 
% outer automorphisms of an infinite  %nilpotent $p$\@group
% %\journal Arch. Math. (Basel) \vol 31 \year 1978 \pag.  %321--328.

[D] \author  M.R. Dixon
\title Automorphisms of nilpotent and related groups, in `Groups '93 
Galway/St. Andrews', vol.\thinspace 1 (eds. C.M.~Campbell et al.), 
London Math. Soc. Lecture
Note Ser. {\bf 211}, Cambridge Univ. Press, 1995, pp.\thinspace 
156--171.

[DG] \author M. Dugas \and R. G\"obel
\title On locally finite $p$\@groups and a problem of Philip Hall's
\journal  J. Algebra\vol 159 \year 
1993 \pag. 115--138.

[Ku] \author D.W. Kueker 
\title Definability, automorphisms, and infinitary languages, in `The 
Syntax and Semantics of Infinitary Languages' (ed.~J.~Barwise), 
Lecture Notes in Math. {\bf 72}, Springer,  Berlin, 1968, 152--165.

[EN] \author P. Hall 
\title `The Edmonton Notes on Nilpotent Groups', Queen Mary College 
Math. Notes, 1969; also in `The collected works  of Philip Hall' 
(eds. K.W.~Gruenberg and J.E. Roseblade), Oxford Science 
Publications, 1988, pp.\thinspace 415--462.

[H] \author J.A. Hulse 
    \title Automorphism towers of polycyclic groups
    \journal J. Algebra \vol  16\year 1970 \pag. 347--398.

[OP] \author O. Puglisi 
\title A note on the automorphism group of a locally finite $p$-group 
\journal Bull. London Math. Soc. \vol 24 \year 1992 \pag. 437--441.

[RobinH] \author D.J.S. Robinson
\title A contribution to the theory of groups with finitely many  
automorphisms
\journal Proc. London Math. Soc. (3) \vol 35 \year 1977 \pag. 34--54. 

[RobinStA]\author D.J.S. Robinson
\title Applications of cohomology to the theory of groups, in 
`Groups---St. Andrews 1981' (eds.~C.M.~Campbell and E.F.~Robertson), 
London Math. Soc. Lecture Note Ser. {\bf 71}, Cambridge Univ. Press, 
1982, pp.\thinspace 46--80.

[Rot] \author J.J. Rotman 
\title `An Introduction to the Theory of Groups' (4th. edn.), 
Springer, Berlin, 1995.

[ST] \author S. Thomas 
\title Complete existentially closed locally finite groups \journal 
Arch. Math. (Basel)  \vol 44 \year 1985 \pag. 97--109.

\ignora{
%[VW] \author M.R. Vaughan-Lee \and J. Wiegold
%\title Countable locally nilpotent groups of finite exponent
%without maximal subgroups
%\journal Bull. London Math. Soc \vol 13 \year 1981 \pag. 45--46.
}

[Z] \author A.E. Zaleski\u\i 
\title An example of torsion-free nilpotent group having no outer  
automorphisms
\journal Mat. Zametki \vol 11 \year 1972  \pag.221--226
\=Math. Notes\vol 11 \year 1972 \pag. 16--19. 

\endgroup
% \par\penalty 10000\vskip2\bigskipamount
\bigskip\bigskip
\noindent 
Authors' addresses:
\par\penalty 10000\vskip\bigskipamount
\line{\frenchspacing\hsize=0.45\hsize\parindent=0pt
\vtop{G. Cutolo\par\vskip 3pt
Universit\`a degli Studi di Napoli ``Federico II",\par
Dipartimento di Matematica e Applicazioni\par
``R. Caccioppoli", \par
Via Cintia --- Monte S.~Angelo,\par
I-80126 Napoli, Italy
\smallskip
{\it e-mail: }cutolo@unina.it}\hss
\vtop{C. Nicotera\par\vskip 3pt
Universit\`a degli Studi di Salerno,\par
Dipartimento di Matematica e Informatica,\par
Via S.~Allende, \par
I-84081 Baronissi (SA), Italy
\bigskip\smallskip
{\it e-mail: }nicotera@matna2.dma.unina.it}}

\bye

%% file: artmacros.tex
\ifx\undefined\artmacrosLetto\let\artmacrosLetto t\else \fi
%
%%%%%% 2. ALTRE MACROS
%
% Aggiungere `%' prima del prossimo rigo, ridefinendo \ToContents
% in qualche modo (ad esempio \let\ToContents\gobbletwo, se non si vuole 
% creare l'indice) in caso sorgano difficolta` con le macros di
% gxeplain.tex.
%
\input geplain

%
%
% entrambi gli stili che seguono sarebbero accettabili:
%
\let\ref\xrefn			% ci si discosta da eplain
\def\rif#1{\unskip~\xrefn{#1}}%
%
% \input geplaltro.tex
%
% Come al solito, meglio poter definire sequenze di controllo `private'
%
\catcode`@\letter
%
%%%%%% 3. FORMATO PAGINA
%
% Formato A4:
\paperheight=297 truemm
\paperwidth=210 truemm
%
% Margini 
% (si potrebbe anche decidere di differenziare i margini destro e sinistro 
% per le pagine pari e dispari)
%
\topmargin = 4 truecm		% <<----
\bottommargin = 4 truecm	% <<----
\leftmargin =28 truemm		% <<----
\rightmargin =28 truemm	% <<----
%
% Posizione del numero di pagina e `running head'.
% Il comando \runninghead={TESTO} specifica il testo da fare apparire,
% non alla prima pagina, in cima al foglio; per default in \smcaps,
% cioe` \caps a 8 pt.
%
\footline={\hss\lower 0.5 truecm \hbox{\sevenrm ---\enskip
		\folio\enskip---}\hss}
\newtoks\runninghead
% \headline={\ifnum\pageno>1
% 		\line{\hss\smcaps\the\runninghead\hss}%
% 			\fi
% 	}%
\headline={\line{\hss
		\ifnum\pageno>1\smcaps\the\runninghead\hss\fi
				}%
	}%

%
%%%%%% 4. SEZIONI
%
\newcount\seznum	% numero della sezione
\font\sezfont=cmbx12 % scaled \magstep2%
%
% Uso: \Sezione {TITOLO}{LABEL}. Potrebbe essere opportuno avere anche
% un parametro per il nome abbreviato in running head
%
%			 TUTTI I PARAMETRI SONO 	<<----
%     				 DA SETTARE		<<----
%
\def\Sezione#1\par{%
	\global\advance\seznum 1	% aggiorna il numero di sezione
	\global\thnum=0 		% azzera quello di enunciato
%	{\let\detto\gobble
%	\edef\t@mp{{#1}}%
%	\xdef\seztit{\expandafter\killtrailingspaces\t@mp}%
%	}%
	   %%%%%%%% TUTTI I FONTS SONO DA SCEGLIERE,
	   %%%%%%%% LE MISURE DA RIGUARDARE
	% deve essere possibile cambiare pagina _prima_ del titolo
	% oppure dopo il primo rigo di testo della nuova sezione,
	% in nessun caso immediatamente dopo il titolo.
	\par\penalty -2000		% invito a cambiare pagina
  			% spazio che precede il titolo     <<----
    \vskip 24 pt plus 12 pt minus 8 pt
		% se non e` possibile il cambio pagina alla precedente 
		% penalty, ne' e` possibile arrivare al primo rigo
		% di testo della nuova sezione, si forza qui un
		% cambio pagina, lasciando in bianco lo spazio 
		% necessario in fondo alla pagina. 
		% Qualora sia invece possibile un cambio pagina
		% ordinario, i comandi nel prossimo rigo
		% non hanno alcun effetto.
	\vfil\penalty 100 \vfilneg
%	\vbox{%
	\let\detto\@sezdetto			% per definire LABEL
	\sezfont %\parindent=2em			% <<---------
%	\item{\the\capnum.\the\seznum.}\ignorespaces#1\par  % titolo + LABEL
\centerline{\the\seznum.\enspace #1\unskip}
%	}%   chiude \vbox
	   % Poiche' \doinpage non funziona all'interno di una \vbox
	   % (perche' genera un \mark che verrebbe altrimenti perso),
	   % i dati per l'indice devono essere spostati fuori dalla \vbox.
%	\ToContents{Sezione}{{\the\capnum.\the\seznum}{\seztit}}%
	\penalty 10000\vskip\bigskipamount			% <<----
	\rm\noindent
	\ignoreblanks			% ignora gli eventuali righi bianchi
}%
\def\@sezdetto#1{\definexref{#1}{\the\seznum}{Sezione}}
%%%%%% 5. ENUNCIATI (numerati e non)
%
%
\newcount\thnum
%
% Uso: \Enunciato{NOME}[(NOTA)][\detto{LABEL}]enunciato\par
% Se il primo carattere (non contando gli spazi) che segue {NOME}
% e` `(', le macros suppongono che a nome segua una NOTA tra parentesi 
% da stampare in \rm, ed anche il punto successivo sara` in \rm.
%
\def\Enunciato#1{%
	\ifdim\lastskip<\bigskipamount\removelastskip\vskip\bigskipamount\fi
	\global\advance\thnum 1
	\begingroup
	\let\detto\@thdetto
	\penalty -500
	\noindent
% 	\bf\the\seznum.\the\thnum.\enspace\ignorespaces#1\unskip
 	\bf#1 \the\seznum.\the\thnum
	\futurenonspacelet\@temp\@enunciato
}%
\def\NoslEnunciato#1{%
	\ifdim\lastskip<\bigskipamount\removelastskip\vskip\bigskipamount\fi
	\global\advance\thnum 1
	\let\detto\@noslthdetto
	\penalty -500
	\noindent
	\bf#1 \the\seznum.\the\thnum
	\futurenonspacelet\@temp\@noslenunciato
}%
% \let\noslhook\begingroup
%
% % \let\@slenunciato\@enunciato
% \let\enfont\sl
% 
%
\def\@enunciato#1\par{%
	\parindent=20pt% 		<<------------  (per \item e \\)
	\def\\##1){\item {$(##1)$}}%	per liste
	\interlinepenalty=200% scoraggia cambi pagina in enunciato   
	% se il primo simbolo incontrato e` `(', allora
	% si assume che l'enunciato sia preceduto da una nota
	% da stampare in \rm.
	\let\@next\ignorespaces
	\ifx\@temp(\space\let\@next\@nota
		\else\unskip.\enskip\fi
	\sl\@next#1\par\endgroup\bigskip
}%
%
% \def@enunciato
%
% \@nota si usa solo nell'ambito di \@enunciato.
%
\def\@nota(#1){\unskip{\rm\ (#1).\enskip}\ignorespaces}%
\def\@noslenunciato{%
	% se il primo simbolo incontrato e` `(', allora
	% si assume che l'enunciato sia preceduto da una nota
	% da stampare in \rm.
	\let\@next\ignorespaces
	\ifx\@temp(\space\let\@next\@nota
		\else\unskip.\enskip\fi
	\rm\@next%
}%
\def\Pf{\ifdim\lastskip=\bigskipamount\removelastskip
	\vskip\medskipamount\fi
	\def\\##1){\par $(##1)$\quad\ignorespaces}%	per liste
	\noindent{\it Proof --- }\rm
}%
\let\Dim\Pf
\def\Th{\Enunciato{Theorem}}%
\def\Prop{\Enunciato{Proposition}}%
\def\Lemma{\Enunciato{Lemma}}%
\def\Cor{\Enunciato{Corollary}}%
\def\Example{\NoslEnunciato{Example}}%
\def\@thdetto#1{\definexref{#1}{\the\seznum.\the\thnum}%
					{Enunciato}\ignorespaces}%
\def\@noslthdetto#1{\definexref{#1}{\the\seznum.\the\thnum}%
					{Enunciato}}%\ignorespaces
%
%\let\detto\undefined
% DA FARE:
%	- \let\lemma\Lemma etc. (?)
%
%%%%%%%	5.	FONTS
%
%
% Uso: \fnt xx XX N n. da`  \font \Nxx=XXn
% Esempio: \fnt cmcsc caps ten 10. equivale a \font\tencaps=cmcsc10
%
\def\fnt#1 #2 #3 #4.{\expandafter\font\csname#3#1\endcsname=#2#4}
%
%
%\font\Bigbf=cmbx12
%\font\Bigit=cmti12
\font\caps=cmcsc10
\font\smcaps=cmcsc10 at 8 pt
%\scriptfont4=\seveni
%\scriptscriptfont4=\fivei
%
% DA FARE:	- \bit per titoli (scegliere la dimensione che serve)
%
%%%%%% 6. BIBLIOGRAFIA
%
\newcount\bibno
%
%% 
 % \ListOfReferences [L1] [L2] .. [Ln] \EndOfReferences
 % definisce, nell'ordine, le labels da L1 a Ln come i
 % numeri da 1 a n allo scopo di riferimenti bibliografici.
 % N.B. Tutto quanto appare nel testo, oltre alle [Li]
 % viene letto da TeX in modo abituale, ma all'interno
 % di un gruppo.
 % Il meccanismo usato comprende:
 % \bibxdef{LABEL}{DEFINITION} (dove LABEL e' `protetta'
 %   per uso nella sola bibliografia e DEFINITION non e'
 %   necessariamente un numero)
 % \bibrefnumber{LABEL}    --> DEFINITION
 % \bibref{LABEL}          --> [DEFINITION]
 %%
\def\bibref#1{[\bibrefnumber{#1}]}%
\begingroup
\catcode`\[\active
\catcode`\@\letter
\gdef\ListOfReferences{\begingroup\catcode`\[\active\relax
	\def[##1]{\global\advance\bibno1\bibxdef{##1}{\number\bibno}}%
	\biblist@
}%
\gdef\bibxdef#1#2{\definexref{@BIBlio:#1}{#2}{biblio}}%
\gdef\bibrefnumber#1{\refn{@BIBlio:#1}}%	
\gdef\biblist@#1\EndOfReferences{#1\endgroup}%
\gdef\references{\parindent=19pt\parskip 2pt plus 1 pt minus 1 pt
		\interlinepenalty 9999\tolerance 1000 \frenchspacing
		\exhyphenpenalty 9999
		\vskip 2.5 \bigskipamount plus \bigskipamount minus \bigskipamount
		\centerline{\sezfont References}
		\penalty 100000\vskip 1,5\bigskipamount
		\bibno=0
	% definizioni	  
	\catcode`\[=13
	\def\and{{\rm and }}%
	\def[##1]{\bibliocheck{##1}\rm\item{\bibref{##1}}}%
	\def\author{\caps\ignorespaces}%
	\def\title{\unskip\rm, \ignorespaces}%
	\def\journal {\unskip, \sl\ignorespaces}%
	\def\vol {\unskip\space\bf\ignorespaces}%
	\def\year {\unskip\space\rm(\ignorespaces}%
	\def\pag.{\unskip) \ignorespaces}%
	\def\={\enskip=\enskip\sl\ignorespaces}% per russi
}%
\gdef\bibliocheck#1{\edef\thel@bel{@BIBlio:#1}%
  % check whether LABEL is already defined.
  \expandafter\ifx\csname\xrlabel{\thel@bel}\endcsname\relax
	% LABEL does not exist
	      \message{^^J>>> ERROR: undefined reference label:
		  `#1'.}
		\else
	% LABEL exists
  	% check whether numbering is correct
  	% (here it is assumed that references are numbered)
		\global\advance\bibno 1
		  \ifnum\bibrefnumber{#1}=\bibno\else%\badnumbermessage
		  \message{^^J>>> ERROR: an inconsistency in
		  the ordering of reference numbers has been found
		  at reference \squarebraket#1] (it was nr. \the\bibno\space at
		  first, now it appears to be nr. \bibrefnumber{#1}).}%
		  \fi
		\fi
}%
\endgroup % ripristinati   \catcode di [ e @
\def\squarebraket{[}% for use in  the message
%
%
%
%%%%%% 7. SIMBOLI E MISCELLANEA
%
% Lettere greche
%
\mathdef\eps{\mathchar"122}%		epsilon tonda (\varepsilon)
\mathdef\phhi{\mathchar"127}%		phi tonda (\varphi)
%
% Altri simboli ordinari
%
%
% simboli di relazione
%
% simboli di operazione binaria
%
% se si vuole sostituire il \wr standard:
%\def\wr{\mathbin{{\rm wr}}}			% intrecciato
%
% operatori 
%
%
% frecce
%
%
% extra
%
%
\catcode`@\other

%% file: geplain.tex
% estende gplain con:
% 	eplain privo di supporto per bibtex
%    +	gxeplain + geplaltro
%    +	formato A4 di default ( 210 x 297 mm )
%
\input gplain

\let\nobibtex t%
\input eplain
\input gxeplain

\input geplaltro

%

%% file: gplain.tex
% vengono definiti:
%
%    +	\italiano, \english
%    +	fonts ams (msam, msbm, eusm), divisi in 
% 	famiglie come da tabella:
%
% 	famiglia		comando		font
% 	____________________________________________
% 	amsafam		amsa		msam
% 	amsbfam=bbbfam	bbb		msbm
% 	eusmfam=tondofam	tondo		eusb
% 	____________________________________________
%    +	\Bbb e \Tondo
%	\fnt e \famiglia 			(da myamsfont.tex)
%    +	diversi simboli,
%	\mathdef, \bbar, \qed e vari caratteri	(da mysymbols.tex)
%    +	\relpenalty=1500 e \binoppenalty=3000
%    +	spaziatura dimezzata per display 
%    +  \ignora
%
%%%%%%%%%%%%%%%%%%%%%%%%%%%%%%%%%%%%%%%%%%%%%
% 
%
%
%
\input myamsfonts

\input mysymbols
\relpenalty=1500           % plain TeX default 500
\binoppenalty=3000         % plain TeX default 700
\abovedisplayskip= 6 pt plus 2 pt minus 3 pt
\belowdisplayskip= \abovedisplayskip
\belowdisplayshortskip= 3 pt plus 1 pt minus 2 pt
\long\def\ignora#1{}

%% file: myamsfonts.tex
% macros per caricare famiglie di fonts e alcuni fonts ams
%
\def\fnt#1 #2 #3 #4.{\expandafter\font\csname#3#1\endcsname=#2#4}%
%
% \famiglia xxx(font) crea la famiglia `xxxfam' 
% con: \textfont\xxxfam=font10, \scriptfont\xxxfam=font7
% e \scriptscriptfont\xxxfam=font5.
%
% per poter definire \famiglia senza aver letto eplain e` necessario
% definire la versione interna di \newfam:
%
\edef\innernewfam{\noexpand\newfam}
\def\famiglia#1(#2){\expandafter\innernewfam\csname#1fam\endcsname
   \fnt #1 #2 ten 10. \fnt #1 #2 seven 7. \fnt #1 #2 five 5.%
   \expandafter\textfont\csname#1fam\endcsname=
		\expandafter\csname ten#1\endcsname
   \expandafter\scriptfont\csname#1fam\endcsname=
		\expandafter\csname seven#1\endcsname
   \expandafter\scriptscriptfont\csname#1fam\endcsname=
		\expandafter\csname five#1\endcsname
}%
%
% famiglia		comando		font
% ____________________________________________
% amsafam		amsa		msam
% amsbfam=bbbfam	bbb		msbm
% eusmfam=tondofam	tondo		eusb
%
\famiglia amsa(msam)%
\famiglia amsb(msbm)%
	\let\bbbfam\amsbfam
	\def\bbb{\tenamsb\fam\bbbfam\relax}%
\famiglia eusm(eusm)%
	\let\tondofam\eusmfam
	\def\tondo{\teneusm\fam\eusmfam\relax}%
%
% Abbreviazioni (e.g., \Bbb Z equivale a {\bbb Z} )
%
\def\Bbb#1{{\bbb #1}}%
\def\Tondo#1{{\tondo #1}}%

%% file: mysymbols.tex
% Alcuni simboli. Presuppone definite:
%
% famiglia		comando		font
% ____________________________________________
% amsafam		amsa		msam
% amsbfam=bbbfam	bbb		msbm
% eusmfam=tondofam	tondo		eusb
%
% (come in myamsfonts.tex)
%
%
% \mathdef\cs{X} definisce \cs come X, con la possibilita` di usarlo sia
% in modo matematico che non. Se in modo non matematico, \cs 
% aggiungera` anche uno spazio. Ad esempio, usi corretti sono:
% 	$3\in\cs$	chiamiamo \cs il...	...detto $\cs$.
% Nell'ultimo caso, `...detto \cs.' non va bene perche' aggiungerebbe
% uno spazio tra X (esito di \cs) e il punto.
%
\def\mathdef#1#2{\def#1{\ifmmode#2\else$#2$ \fi}}
%
% Per definire \ordchardef, \binchardef  etc... c'e` bisogno
% di convertire numeri in forma decimale nella forma esadecimale.
% La difficolta` e` data dalla necessita` di uso con \edef.
%
\catcode`@=11
%
% \hexdigit <num> esprime <num> in una cifra esadecimale, se 
% 0 =< <NUM> < 16. La definizione e` cosi` stupida per poter funzionare
% sotto \edef
%
\def\hexdigit#1{\ifnum#1<10 \number#1\else
	 \ifnum#1=10 A\else
	  \ifnum#1=11 B\else
	   \ifnum#1=12 C\else
	    \ifnum#1=13 D\else
	     \ifnum#1=14 E\else
	      \ifnum#1=15 F\else
		\error{out of range for \hexdigit}%
	\fi\fi\fi\fi\fi\fi\fi}%
%
% Quasi come sopra, per numeri a due cifre esadecimali.
% \hextwodigits NUM TOKS definisce il token register TOKS come l'espressione
% esadecimale del numero NUM a due cifre esadecimali.  Il risultato in TOKS
% e` sempre in due cifre (ad es. `0A', non `A'), perche' questo e` il
% formato che serve nelle macros successive.
%
\def\hextwodigits#1#2{\count255=#1 \chardef\t@mpnum#1%
		\divide\count255 by 16
		\edef\t@mp{\noexpand#2{\hexdigit{\count255}}}\t@mp
		\multiply\count255 by -16
		\advance\count255 by \t@mpnum
		\edef\t@mp{\noexpand#2{\the#2 \hexdigit{\count255}}}%
		\t@mp
}%
%
% \defchardef n{tipo} definisce \tipochardef, che ha il seguente uso: 
% \tipochardef\cs=\xxxfamN	equivale a 
% \mathchardef\cs="nFXX		dove F e` la cifra esadecimale 
% corrispondente alla famiglia \xxxfam e XX e` la forma esadecimale
% (in due cifre, vedi sopra) del numero (decimale) N.
%
\def\defchardef#1#2{%
	\expandafter\edef\csname#2chardef\endcsname##1=##2{%
		\edef\noexpand\t@mpdef{%
			\mathchardef\noexpand\noexpand ##1="#1%
		\noexpand\hexdigit{##2}}\afterassignment
		\noexpand\@finechardef\noexpand\count255=}%
}%
\def\@finechardef{\hextwodigits{\count255}{\toks0}%
		\expandafter\t@mpdef\the\toks0\relax
}%
\catcode`@=12
\defchardef0{ord}%	definisce \ordchardef
\defchardef1{Op}%	definisce \Opchardef
\defchardef2{bin}%	definisce \binchardef
\defchardef3{rel}%	definisce \relchardef
%
%
% Simboli in blackboard bord
%
\mathdef\Z{{\Bbb  Z}}%		interi
\mathdef\N{{\Bbb  N}}%		
\mathdef\No{{\Bbb  N_0}}%	naturali
\mathdef\Q{{\Bbb  Q}}%		razionali
\mathdef\S{{\Bbb S}}%		S per gruppi simmetrici
%
% Altri simboli ordinari
%
\ordchardef\vuoto=\bbbfam63		% insieme vuoto
\mathdef\phhi{\varphi}%
\mathdef\eps{\varepsilon}%
%
%
% simboli di operazione binaria
%
\def\setminus{\mathbin{\raise 1pt\hbox{\bbb\char114}}}%%% setminus
\binchardef\semid=\bbbfam111				% semidiretto
\let\semidirect\semid					%    ,,
%
% simboli di relazione
%
\relchardef\Re=\tondofam82		% R di relazione
\relchardef\n=\amsafam67		% normale
\relchardef\notn=\bbbfam54		% non normale
			% subnormale
\let\iso\simeq				% isomorfo
\let\st\mid				% tale che
%
% operatori 
%
\def\gen#1{\langle#1\rangle}%			Sottogruppo generato da #1
%		ordine di #1
\def\im{\mathop{\rm im}\nolimits}
%
% frecce
%
\def\mono{\hskip2pt\hbox{{\def\joinrel{\mathrel {\mkern
    -4mu}}$\raise1,1pt\hbox{$\scriptscriptstyle
     >$}\kern-5pt\longrightarrow$}}} 
\def\epi{\hskip2pt\hbox{\rlap{$\longrightarrow$}}\kern-1pt
    \longrightarrow}

%
% extra
%
\def\@{\/\hbox{-}}%			trattino che non va a capo
\def\({\bigl(}%
\def\){\bigr)}%				parentesi tonde grandi
\def\[{\bigl[}%
\def\]{\bigr]}%				parentesi quadre grandi
% \bbar X crea una barretta piu` lunga di quella generata da \bar su X. 
% Da usare in luogo di \bar per le maiuscole.
\def\bbar#1{\kern 2pt\rlap{$\bar{\phantom{#1}}$}\kern-2pt\bar#1}
% \qed e` una simbolo intelligente di fine dimostrazione.
% Inserisce un quadratino vuoto (\qedbox) alla fine del rigo. Se lo 
% spazio non e` sufficiente, puo` creare un rigo bianco ed inserire 
% \qedbox alla sua fine, ma cerchera` prima di ridisegnare il paragrafo
% affinche' cio` non accada. Volendo si puo` ridefinire \qedbox.
\def\qedbox{$\rlap{$\sqcap$}\sqcup$}%
\def\qed{\nobreak\hfill\penalty500 \hbox{}\nobreak\hfill\qedbox\par\medbreak}%

%% file: gxeplain.tex
% modifiche (aggiunte) a:
% xeplain.tex: macros for nonformatting.  Written 1989--94 by (mostly)
% Karl Berry.  These macros are in the public domain.
%
%
%\let\nobibtex=t
%\input eplain
\catcode`@\letter
\let\@plainwlog\wlog
\let\wlog\gobble
%
%
%%%%%%%%		Deferred commands
%
% The scope of the macros in this section is to make possible typing
% \doinpage{COMMANDS} to have the COMMANDS read and executed by Tex only
% during the next output routine. The introduction of \everyoutput could
% be of other use as well.
% The idea is to use \mark to keep track of the number of commands issued
% in the current page. This has the disadvantage that \mark should not be 
% used in other ways. Therefore it is disabled here. (It seems though that
% whatever you wanted to do with \mark can still be done with \doinpage,
% possibly less efficently.)
%
% The real problem with \doinpage is that it could fail to work in the
% cases when \mark goes lost, for instance in internal vertical mode
% (see the TeXbook, p. 259). Hence \doinpage should not be used inside
% complex box structures.
% However, since the only relevant thing to the output routine is what
% the last mark issued in the page is, every `well placed' \doinpage
% command will rescue all misplaced one preceding it in the same page.
% For instance, \vbox{...\doinpage{A}...\doinpage{B}...
%		      ...\doinpage{W}...}\penalty10000\doinpage{Z}
% is perfectly safe.
%
\newtoks\everyoutput       % read at every page output
\newcount\@todopg          % number of deferred commands to be read
\newcount\@pgdone          % number of deferred command read
\newtoks\dopgtoks          % list of deferred commands to be read
%
% disable mark to avoid conflicting usage.
%
\let\plainmark\mark
\let\mark\undefined
%
% add \everyoutput to the output routine. The default value is \doinpageloop
%
\output={\the\everyoutput\plainoutput}%
\everyoutput={\doinpageloop}%
%
% \doinpage{ARG} normally causes the token list ARG to be read at the next
% end of page by \doinpageloop. The counter \@todopg is advanced by one 
% and marked. In this way, a list of commands is produced, and \botmark
% contains a number corresponding to the most recent command issued in the
% current page.
%
\def\doinpage#1{%
   \global\advance\@todopg\@ne
   \plainmark{\the\@todopg}%
      % appends #1 followed by a separator to \dopgtoks
   \global\dopgtoks\expandafter{\the\dopgtoks#1@@DoNext@@}%
}%
%
%
% \doinpageloop make TeX read the commands issued by \doinpage in the 
% current page (i.e., up to the one corresponding to \botmark).
% Each command is expunged from the list \dopgtoks after it is read.
% (The counters \@pgdone and \@todopg are different to this respect: while
% \dopgtoks contains only the commands yet to be read, both counters contain
% the number of all deferred commands read or to be read in the document.)
%
\def\doinpageloop{%
   \loop
   \ifnum 0\botmark>\@pgdone
   % the `0' before \botmark lets \ifnum read a number even if \botmark
   % is missing.
     \global\advance\@pgdone\@ne
     \doOnepg
   \repeat
}%
%
% \doOnepg reads the first command (=argument of \doinpage) from \dopgtoks
% and discards it from the latter.
%
\def\doOnepg{\expandafter\@dopgpickfirst\the\dopgtoks @@PGTOKS_STOP@@}
\def\@dopgpickfirst#1@@DoNext@@#2@@PGTOKS_STOP@@{\global\dopgtoks{#2}#1}
%
%
%
%
%
%%%%%%%%%	  Cross-references.
%
% Definition and information about a reference labelled LABEL
% are stored in the control sequence LABEL_x and its properties `class' 
% and `status'.
% 
% If the document contains a reference to LABEL preceding its definition,
% the latter is recorded in \jobname.aux, called \auxfile in the macro.
% 
% If LABEL_x is defined, its status is one of the following:
% 
% d	(definitive) if LABEL_x has already been defined in the document
%	(as opposite to \auxfile);
% a	(aux) if LABEL_x has not yet been defined in the document, and no
%	reference to LABEL has yet been made (so LABEL has only been found
% 	by TeX in \auxfile).
% r	(referred to) if LABEL_x has not yet been defined in the document
%	but a reference to LABEL has already been made. In this case
%	LABEL_x either has been defined in \auxfile or has a dummy 
%	definition.
%	
% Using `status' allows to avoid running eplain twice to get 
% cross-references right, if possible, getting a warning otherwise.
% A warning is also printed if a label is defined twice. Moreover,
% \auxfile is created only if really necessary (though this feature can de 
% disabled with \wantauxfile).
% 
%
\newif\if@wantauxfile
\def\wantauxfile{\immediate\openout\auxfile=\jobname.aux%
		\global\@wantauxfiletrue}%
%
% When a label isn't defined, we only want to complain if
% \xrefwarningtrue; btxmac uses \if@citewarning for this, so we have to
% reuse that name.  We can't just say \let\ifxrefwarning =
% \if@citewarning, since then changes to the latter won't be reflected
% in the former.  On the other hand, we have to have a true \if...
% command, so \if's and \fi's match properly.  What a mess.
%
\let\ifxrefwarning = \iftrue
\def\xrefwarningtrue{\@citewarningtrue \let\ifxrefwarning = \iftrue}%
\def\xrefwarningfalse{\@citewarningfalse \let\ifxrefwarning = \iffalse}%
%
% If a second run is necessary to get cross-references right, a warning
% is printed, but no more than once per run.
%
\newif\if@rerunyelled
\def\@rerunmessage{\if@rerunyelled\else\message{%
	>>> YOU MUST RUN TEX AGAIN TO GET CROSS-REFERENCES RIGHT^^J}%
	\global\@rerunyelledtrue
	\fi}%
\def\@twicedefinedmessage #1{\message{%
	line \linenumber - label `#1' already defined.
	New definition ignored^^J}}%
% 
%
% \definexref{LABEL}{DEFINITION}{CLASS} defines a cross-reference named
% LABEL of label class CLASS to be DEFINITION.  (Or LABEL can be a
% control sequence; it's expanded to get the label text.)  To get a
% possible page number right, we have to write the definition out to the
% auxiliary file, instead of only defining it directly.
%
% Recall that \xrlabel{LABEL} expands to a cross-reference internal name,
% that is, LABEL_x.
%
\def\definexref#1#2#3{%
  % Remember what we're given; it might be `\@optionalarg' which
  % \readauxfile trashes. Moreover, LABEL could contain control sequences
  % like \relax to be expanded before \csname can complain about them.
  % (No loss of generality here since \csname will fully expand the label 
  % anyway)
  \edef\thel@bel{#1}%
  %
  % check whether LABEL (#1=\thel@bel) is already defined. If not, define 
  % it with status `d'. The status will be `d' in all cases.
  %
  \expandafter\ifx\csname\xrlabel{\thel@bel}\endcsname\relax
    \@definelabel{\thel@bel}{#2}{#3}{d}%
	\if@wantauxfile
		\@writedeftoaux{\thel@bel}{#2}{#3}%
	\fi
  \else
  % take the status of LABEL_x and consider the different possibilities.
    \edef\@labelstatus{\getproperty{\xrlabel{\thel@bel}}{status}}%
    %
    % if LABEL_x has been defined only in \auxfile and never referred to
    % in the document so far, redefine it (and don't write it to \auxfile,
    % unless \wantauxfile is in force).
    %
    \if \@labelstatus a%
	\@definelabel{\thel@bel}{#2}{#3}{d}%
	\if@wantauxfile
		\@writedeftoaux{\thel@bel}{#2}{#3}%
	\fi
    \else
    %
    % if LABEL_x has already been defined in the document, print a warning
    % (unless \xrefwarning is false) and do nothing else. (One could check
    % whether the `old' and the `new' definitions agree, but what for?)
    %
	\if\@labelstatus d%
		\ifxrefwarning\@twicedefinedmessage{#1}\fi
	\else
	%
	% LABEL is defined here for the first time in the document, but
	% has been already referred to. Redefine it and write the 
	% definition to \auxfile,
	% and check whether the new definition agrees with the previous.
	% If not, a second run of the document is needed, and the
	% suitable warning is printed.
	% However, this check takes place only if such a warning has not
	% yet been printed.
	% Finally, since \auxfile has been opened, there is no point in
	% avoiding writing definitions to it, and \wantauxfile is set.
	%
	   \if@wantauxfile\else\wantauxfile\fi
	   \@writedeftoaux{\thel@bel}{#2}{#3}%
	   \if@rerunyelled
	   \else
		\TestLabel{\thel@bel}{#2}{#3}%
		\ifx\TestLabelResult f% 
			\@rerunmessage
		\fi
	   \fi
	   \@definelabel{\thel@bel}{#2}{#3}{d}%
	\fi
    \fi
  \fi
}%
%
% In \auxfile definitions are given by \auxdefinexref, with similar
% parameters for LABEL, DEFINITION and CLASS. If LABEL_x was already
% defined, the status is not changed, otherwise it is set to a.
% If the status is d, the new definition is is ignored.
%
\def\auxdefinexref#1#2#3{%
  \edef\thel@bel{#1}%
  \if\getproperty{\xrlabel{\thel@bel}}{status}d%
  \else
    \if\getproperty{\xrlabel{\thel@bel}}{status}r\relax
	% the \relax here is useful to avoid \if expanding too far:
	% if status of #1 is not defined then \if compares r and the
	% following token, after expansion. \relax is not needed after
	% d in the previous \if  condition because d is directly followed
	% by \else, so that if the `status' token is missing \if is happy
	% with the condition being false.
	\@definelabel{\thel@bel}{#2}{#3}{r}%
    \else
	\@definelabel{\thel@bel}{#2}{#3}{a}%
    \fi
  \fi
}%
%
% \@definelabel{LABEL}{DEFINITION}{CLASS}{STATUS} actually defines LABEL 
% of label class CLASS and status STATUS to be DEFINITION.
%
\def\@definelabel#1#2#3#4{%
  % Define the control sequence.
  \expandafter\xdef\csname\xrlabel{#1}\endcsname{#2}%
  %
  % Remember what kind of label this is.
  \global\setproperty{\xrlabel{#1}}{class}{#3}%
  \global\setproperty{\xrlabel{#1}}{status}{#4}%
}%
%
% \@writedeftoaux#1#2#3 writes `\auxdefinexref{#1}{#2}{#3}' to \auxfile
%
\def\@writedeftoaux#1#2#3{%
  % Remember what we're given; it might be `\@optionalarg', which
  % \readauxfile trashes.  (No loss of generality here, since \csname
  % will fully expand the label anyway.)
  \edef\t@mp{#1}%
  %
  % Be sure we've read the aux file before we zap it:
  \readauxfile
  %
  % When we read in the aux file next time, define the label:
  \edef\@wr{\noexpand\immediate\noexpand
	    \writeaux{\string\auxdefinexref{\t@mp}{#2}{#3}}}%
  \@wr
  \ignorespaces
}%
%
% \TestLabel{LABEL}{DEFINITION}{CLASS} tests whether DEFINITION and CLASS
% respectively are definition and class of the label LABEL. The result
% (t/f) is stored in \TestLabelResult.
%
\def\TestLabel#1#2#3{%
 % \t@mp is set to D@@@C where D is the definition and C is the class of #1
 \edef\thel@bel{#1}%
 \edef\t@mp{\csname\xrlabel{\thel@bel}\endcsname 
	@@@\getproperty{\xrlabel{\thel@bel}}{class}}%
 \edef\t@@mp{#2@@@#3}%
 \ifx\t@mp\t@@mp
	\let\TestLabelResult t%
 \else
	\let\TestLabelResult f%
 \fi
}%
%
% \xrefn{LABEL} typesets a reference to the label LABEL, keeping
% track of the needed `status' information.
%
\def\xrefn#1{%
  % if LABEL_x is not defined check first \auxfile. Then, if LABEL_x is
  % still undefined, give it a dummy definition with status `r'
  % and, unless disabled, print the corresponding warning.
  %
  \expandafter \ifx\csname\xrlabel{#1}\endcsname\relax
    \readauxfile
    \expandafter \ifx\csname\xrlabel{#1}\endcsname\relax
      \if@citewarning
	\message{\linenumber Undefined label `#1'.^^J}%
      \fi
      \expandafter\gdef\csname\xrlabel{#1}\endcsname{%
	`{\tt\escapechar = -1
	  \expandafter\string\csname#1\endcsname}'%
      }%
    \fi
    \global\setproperty{\xrlabel{#1}}{status}{r}%
  \fi
  \csname\xrlabel{#1}\endcsname % Always produce something.
  \if\getproperty{\xrlabel{#1}}{status}a\relax
	\setproperty{\xrlabel{#1}}{status}{r}%
  \fi
}%
%
% \refn is just a synonym.
%
\let\refn = \xrefn
%
% \readauxfile looks for the \auxfile. If it finds it, then it reads it.
% this will done at most once in the run, thanks to the \if@auxfiledone
% check.
%
\def\readauxfile{%
	\if@auxfiledone\else
	\global\@auxfiledonetrue
	\@testfileexistence{aux}%
		\if@fileexists
			\begingroup\endlinechar=-1\catcode`@=11
			\input \jobname.aux \endgroup
		\fi
	\fi
}%
%
%
% \xrdef{LABEL} defines the label LABEL as the current page number.
%
%
%
% The eplain control sequence \@eqdefn, for defining equation
% number references, uses \@definelabel. Eplain \@definelabel requires
% three paramenters, whereas our \@definelabel requires four. So we
% must re-define \@eqdefn. It is very straightforward.
%
\def\@eqdefn#1#2{%
	\definexref{#1}{#2}{eq}%
}%
%
%%%%%%%%%	Table of Contents
%
% Data for Table of Contents are stored in the token list \ToCtoks 
% and normally read by \ToC.
%
\newtoks\ToCtoks
%
% \ToContents{ENTRY_TYPE}{STUFF} appends to \ToCtoks the tokens
% \ToCENTRY_TYPE{pag}STUFF, where {pag} is the page number.
% The insertion of the page number makes a call to \doinpage necessary.
% Both ENTRY_TYPE and STUFF are expanded before being appended to 
% \ToCtoks.
% Under normal usage, \ToCtoks will be read by \ToC.
%
% In order to make this work, for every ENTRY_TYPE appearing in the document,
% the control sequence \ToCENRTY_TYPE must be defined when \ToC is 
% invoked. Moreover, \ToCENRTY_TYPE must require parameter text
% of the form #1_OTHER_ where #1 will stand for the page number and 
% every occurence of STUFF corresponds to _OTHER_ (not counting the
% braces surrounding STUFF). 
% E.g., \ToContents{Chapter}{{num}{A title}} will generate the entry
% \ToCChapter{P}{N}{A title}, where P and N are the page and the chapter
% numbers, hence the matching definition will have the form
% \def\ToCChapter#1#2#3{...}.
%
\def\ToContents#1#2{%
		\doinpage{\edef\t@mp{%
			\global\noexpand\ToCtoks\noexpand\expandafter
			{\noexpand\the\ToCtoks\expandonce\csname 
			  ToC#1\endcsname{\folio}#2}%
			  }%
		\t@mp}%
}%
%
% \ToCformat is an example format for the ToC page(s). Every entry makes
% a paragraph. Page numbers are printed flush right, preceded by leading 
% dots. Lines without page numbers are printed ragged-right, of length
% between 0.6 and 0.9 times \hsize. Page number is set to -1 (which
% will be printed as `i').
%
\def\ToCformat{%
	\rightskip = .25\hsize plus .15\hsize minus .15\hsize
	\parfillskip=-\rightskip
}%
%
% \ToC is defined just to have one such control sequence working. Each
% document using Table of Contents will probably need to redefine \ToC
% to add some more specific feature.
%
% \ToCheading (the heading to appear on top of the first ToC page) is 
% defined apart to make easier to redefine it for non-english documents.
%
\def\ToCheading{\bf Table of Contents}
\def\ToC{%
	\vfill\eject\centerline{\ToCheading}\pageno=-1
	\vskip 3\bigskipamount plus \bigskipamount minus \bigskipamount
	\ToCformat
	\the\ToCtoks
}%
\let\wlog\@plainwlog
\catcode`@\other

%% file: geplaltro.tex
\catcode`@\letter
%
%%%%%%%		while-false loops
%
% Plain TeX \loop ... \repeat command does not allow constructions like
% \loop A \if.. \else B\repeat. The command \whilefalse\if..\do B \repeat
% allows something close. It will execute B until the condition in \if..
% turns out to be true. 
% The definition mimicks that of plain-TeX \loop.
% #1 is a \if.. condition, #2 is the command to be repeatly executed.
%
\def\whilefalse#1\do#2\repeat{\def\condition{#1}\def\body{#2}\wfiterate}
\def\wfiterate{\condition\let\wfiterate\relax\else\body\fi\wfiterate}
%
%
%%%%%%%		\readfile
%
% \readfile FILE in TOKS stores the content of the file FILE in the 
% \toks register TOKS. Of course end-of-line characters are translated
% into spaces and empty lines into \par's. A \par is added as last line.
%
% WARNINGS:
% - Since the routine uses \read, it will not work if FILE contains an 
%   outer macro.
% - never let \toks0 as #2 (TOKS)
%
\def\readfile #1 to #2{%
	\immediate\openin0=#1
	% loops over \ifeof0 being false
	\whilefalse\ifeof0\do
	\read0 to \@temp
	\toks0=\expandafter{\@temp}%
	#2=\expandafter\expandafter\expandafter
		{\expandafter\the\expandafter#2\the\toks0}%
	\repeat
}%
%
%
%%%%%%%		Bibliography
%
% Simple macros for getting bibliographical reference numbers for
% non-bibtex-users.
% Using the cross-reference mechanism of e-plain obviously works, but
% certainly causes the document to be read twice. The following macros
% are designed to make TeX learn the correct numbering of reference items
% just at the beginning of the job.
% 
\catcode`@\letter
\let\@plainwlog\wlog
\let\wlog\gobble
%
%
% \bibliod@ne is true if the list of the bibliography labels has already
% been read (at most one such list may be read in a TeX job).
% \bibn@ is a counter used to assign conecutive reference numbers to the
% labels.
%
\newif\ifbibliod@ne
\newcount\bibn@
%
%
% Bibliograpic references of labels (or biblabels), L1, L2, etc. 
% are defined as ordinary eplain cross references of labels L1@bib, 
% L2@bib, etc. They may be assigned in two different ways:
%
% 1.
% \definebibrefs text [L1] text [L2] ... text [Ln] ... \endbibrefs 
% defines L1, L2, ... Ln as biblabels of values succesive integers
% starting with 1.
% All text non contained between brackets is ignored.
%
% 2.
% Alternatively, one can type \bibliofile FILE (without \endbibnumbers).
% This makes TeX input the file FILE, where presumably the list of 
% references is stored, and read its content as if it where between a
% \definebibrefs ... \endbibrefs pair.
% 
% In both cases the bracketed items [L1], [L2] etc must be sorted 
% manually by the user. 
%
%
% \definebibrefs starts a loop in which \@definebibrefs describes the
% action to be executed at each pass. The loop ends when \@definebibrefs
% finds `[_]' as argument, where the subscript has category 3 (!). This,
% of course, just to keep the switch private.
%
\begingroup
\catcode`\_=3
\xdef\funny@sb{_}
\global\long\def\definebibrefs #1\endbibrefs{%
	% do not define any biblabel, if this has already been done.
	\ifbibliod@ne
		\@bibtwicemessage
	\else
	% start the loop after appending the ending item (`[_]') to the
	% list of biblabels to be defined.
		\@definebibrefs#1[_]%
	\fi
}%
\endgroup
\long\def\@definebibrefs #1[#2]{%
	\global\bibliod@netrue
	% #1 is discarded. 
	% Check whether #2 is the ending `_'...
	\edef\t@mp{#2}%
	\ifx\t@mp\funny@sb
	\else
		% ... if not define the biblabel #2, of class bib
		\global\advance\bibn@\@ne
		\definexref{#2@bib}{\number\bibn@}{bib}%
		% proceed with the loop. \expandafter hides the next \fi
		% to the argument of \@definebibrefs.
		\expandafter\@definebibrefs
	\fi
}%
\def\bibliofile#1 {%
	% keep in mind the file with the bibliography, for later use
	\gdef\@bibliofile{#1}%
	\ifbibliod@ne
		\@bibtwicemessage
	\else
		\readfile #1 to {\toks1}
		\expandafter\definebibrefs\the\toks1 \endbibrefs
	\fi
}%
\def\brefnum#1{\refn{#1@bib}}%
\def\bref[#1]{[\brefnum{#1}]}%
%
%
%
%
%%%%%%%%%		Ignoring blanks
%
% after \ignoreblanks TeX skips all spaces and empty lines (and explicit
% \par tokens).
%
\def\ignoreblanks{\futurenonspacelet\t@mp\@ignblks}%
%
% Parameter #1 in the definition of \g@on just in case a space has
% somehow appeared in between.
%
\def\@ignblks{\ifx\t@mp\par
	\def\g@on##1\par{\ignoreblanks}\else\let\g@on\relax
	\fi
	\g@on
}%
%
% If TOKS is a list of tokens consisting in TOKS1 followed by (ordinary)
% spaces, in such a way that TOKS1 does not end with a space, then
% \killtrailingspaces{TOKS} produces TOKS1. A proviso on TOKS1 is that
% it must not contain 5 consecutive space tokens (which does not seem to
% be particularly likely anyway).
%
% The trick is that ` \t@mp1 ' produces two consecutive spaces after 
% expansion. \@killsp will swallow everything from the first occurrence
% of 5 consecutive spaces; \killtrailingspaces puts 5 of them just after 
% TOKS -to make sure that \@killsp will find at least them- before
% calling \@killsp. The `@' character delimits the effect of \@killsp.
%
\def\t@mp1{}
\edef\@temp{\def\noexpand\@killsp##1 \t@mp1 \t@mp1 \t@mp1 \t@mp1 ##2@{##1}%
	\def\noexpand\killtrailingspaces ##1{%
	\noexpand\@killsp##1 \t@mp1 \t@mp1 \t@mp1 \t@mp1 @}}
\@temp
%
%
%%%%%%%	NOLIST
%
% sort of like \noalign for alignments, \nolist{STUFF} inside a eplain 
% list puts STUFF in between the list as `non listed material',
% that is without indentation (or with the indentation of the previous level 
% of listing, in case of nesting) and with proper vertical spacing.
% STUFF may consist of several paragraphs and is put inside a group.
%
\long\def\nolist#1{%
	\par\vskip\interitemskipamount
	\begingroup
	\advance\leftskip by -\listleftindent
	\advance\leftskip by -\parindent
	\advance\rightskip by -\listrightindent
	\noindent
	\ignorespaces#1\par\endgroup
}%
\let\wlog\@plainwlog
\catcode`@\other